\documentclass[english,12pt]{smfart}
\usepackage{amsfonts}
\usepackage{amscd}
\usepackage{graphicx}
\usepackage{calrsfs}
\usepackage{yfonts}
\usepackage[all]{xypic}
\usepackage[all]{xy}
\usepackage[francais, english]{babel}
\usepackage[latin1]{inputenc}

\CompileMatrices

\swapnumbers

\def\build#1_#2^#3{\mathrel{
\mathop{\kern 0pt#1}\limits_{#2}^{#3}}}

\let\cal\mathcal

\def\11{{\mathbf 1}}
\def\AA{{\mathbf A}}

\def\CC{{\mathbf C}}

\def\PP{{\mathbf P}}
\def\QQ{{\mathbf Q}}
\def\RR{{\mathbf R}}
\def\SS{{\mathbf S}}

\def\ZZ{{\mathbf Z}}

\def\Lambdal{\Lambda^{\ell oc}}
\def\Var{{\hbox{Var}}}

\mathchardef\alphag="7C0B \mathchardef\betag="7C0C
\mathchardef\gammag="7C0D \mathchardef\deltag="7C0E
\mathchardef\varepsilong="7C22 \mathchardef\varphig="7C27
\mathchardef\psig="7C20 \mathchardef\zetag="7C10
\mathchardef\epsilong="7C0F \mathchardef\rhog="7C1A
\mathchardef\taug="7C1C \mathchardef\upsilong="7C1D
\mathchardef\iotag="7C13 \mathchardef\thetag="7C12
\mathchardef\pig="7C19 \mathchardef\sigmag="7C1B
\mathchardef\etag="7C11 \mathchardef\omegag="7C21
\mathchardef\kappag="7C14 \mathchardef\lambdag="7C15
\mathchardef\mug="7C16 \mathchardef\xig="7C18
\mathchardef\chig="7C1F \mathchardef\nug="7C17
\mathchardef\varthetag="7C23 \mathchardef\varpig="7C24
\mathchardef\varrhog="7C25 \mathchardef\varsigmag="7C26
\mathchardef\Omegag="7C0A \mathchardef\Thetag="7C02
\mathchardef\Sigmag="7C06 \mathchardef\Deltag="7C01
\mathchardef\Phig="7C08 \mathchardef\Gammag="7C00
\mathchardef\Psig="7C09 \mathchardef\Lambdag="7C03
\mathchardef\Xig="7C04 \mathchardef\Pig="7C05
\mathchardef\Upsilong="7C07

\newtheorem{ex}[subsubsection]{Example}
\newtheorem{thm}[subsubsection]{Theorem}
\newtheorem{theorem}[subsection]{Theorem}
\newtheorem{lem}[subsubsection]{Lemma}
\newtheorem{lemma}[subsubsection]{Lemma}
\newtheorem{cor}[subsection]{Corollary}
\newtheorem{prop}[subsubsection]{Proposition}

\newtheorem{CCF}[subsection]{Cauchy-Crofton formula}
\newtheorem{LCCF}[subsubsection]{Local Cauchy-Crofton formula}
\newtheorem{MLCCF}[subsubsection]{Multidimensional local Cauchy-Crofton formula}

\theoremstyle{definition}
\newtheorem{definition}[subsubsection]{Definition}
\newtheorem{defi}[subsection]{Definition}
\newtheorem{example}[subsection]{Example}

\newtheorem{def-prop}[subsection]{Proposition-Definition}
\newtheorem{def-theorem}[subsection]{Theorem-Definition}
\newtheorem{def-lem}[subsection]{Lemma-Definition}

\theoremstyle{remark}
\newtheorem{remark}[subsubsection]{Remark}
\newtheorem{rem}[subsection]{Remark}
\newtheorem{rems}[subsection]{Remarks}
\newtheorem{remabout}[subsection]{Remarks about $(1')$ and $(2')$}
\newtheorem{remarks}[subsubsection]{Remarks}
\newtheorem{r and comp}[subsection]{Remarks and computations}
\newtheorem{notation}[subsubsection]{Notation}
\newtheorem{nota}[subsection]{Notation}
\newtheorem{question}[subsubsection]{Question}
\newtheorem{questions}[subsubsection]{Questions}

\theoremstyle{plain}

\numberwithin{equation}{subsection}

\def\boxit#1#2{\setbox1=\hbox{\kern#1{#2}\kern#1}%
\dimen1=\ht1 \advance\dimen1 by #1 \dimen2=\dp1 \advance\dimen2 by
#1
\setbox1=\hbox{\vrule height\dimen1 depth\dimen2\box1\vrule}%
\setbox1=\vbox{\hrule\box1\hrule}%
\advance\dimen1 by .4pt \ht1=\dimen1 \advance\dimen2 by .4pt
\dp1=\dimen2 \box1\relax}

\renewcommand{\theequation}{\thesubsection.\arabic{equation}}

\DeclareMathOperator{\jac}{jac}

\def\R{{ \mathbb R}}
\def\Z{{ \mathbb Z}}
\def\L{{ \mathbb L}}
\def\C{{ \mathbb C}}
\def\P{{ \mathbb P}}
\def\N{{ \mathbb N}}

\newcommand{\eps}{\varepsilon}

\date{\today}

\begin{document}

\title[\tiny{Deformation of singularities and additive invariants}]{\rm Deformation of singularities and additive invariants}{}
\author{Georges COMTE}
\address{Laboratoire de Math\'ematiques de l'Universit\'e de Savoie, UMR CNRS 5127,
B\^atiment Chablais, Campus scientifique, 
73376 Le Bourget-du-Lac cedex, France}
\email{georges.comte@univ-savoie.fr}
\urladdr{http://gc83.perso.sfr.fr/}

\begin{abstract} In this survey on local additive invariants of real and complex definable singular germs we systematically 
present classical or more recent invariants of different nature as emerging from a tame degeneracy principle. For this goal, 
we associate to a given 
singular germ a specific deformation family whose geometry degenerates in such a way that it eventually gives rise to a list 
of invariants attached to this germ. 
Complex analytic invariants, real curvature invariants and motivic type invariants are encompassed under this point of view.  
We then explain how all these invariants are related to each other as well as we propose a general conjectural principle 
explaining why such invariants have to be related. 
This last principle may appear as the incarnation in definable geometry of deep finiteness results of convex geometry, according 
to which additive invariants in convex geometry are very few.  
   
\end{abstract}

\maketitle

\section*{Introduction}
A beautiful and fruitful principle occurring in several branches of mathematics consists in deforming the object under consideration in order to let appear some invariants attached to this object. In this deformation process, the object $X_0$ to study becomes the special fibre of a deformation family $(X_\eps)_\eps$ where each fibre $X_\eps$  
approximates $X_0$, from a topological, metric or geometric point of view, depending on the nature of the invariant that one aims for $X_0$ through $(X_\eps)_\eps$. For instance in Morse theory, where a smooth real valued function $f:M\to \RR$ on a smooth manifold $M$ is given, such that $f^{-1}([f(m)-\eta, f(m)+\eta]$ contains no critical point of $f$ but $m$, the homotopy type of $f^{-1}(]-\infty, f(m)+\eta])$ is given by the homotopy type of 
$f^{-1}(]-\infty, f(m)-\eps])$, for any $\eps$ with $0<\eps\le 
\eta$, plus a discrete invariant attached to $f$ at $m$, namely the index of $f$ at $m$. In this case the family $(f^{-1}(]-\infty, f(m)-\eps]))_{0<\eps\le \eta}$ has the same fibres, from the differential point of view, and approximates the special set
$f^{-1}(]-\infty, f(m)+\eta])$, from the homotopy type point of view, up to some additional discrete topological invariant. 
Another instance of this deformation principle can be found in tropical geometry, where a patchwork polynomial embeds a complex curve $X_0$ of the complex torus $(\CC^*)^2$ into a family $(X_\eps)$ of complex curves. This family may be viewed as a curve $\cal X$ on the non archimedean valued field $\CC((\eps^\RR))$ of Laurent series with exponents in $\RR$. Then, by a result of Mikhalkin and Rullg\aa rd (\cite{Mik} and \cite{Rul}), the amoebas family ${\cal A}(X_\eps)$ of $(X_\eps)$ has limit (in the Hausdorff metric)
the non archimedean amoeba ${\cal A}(\cal X)$ of $\cal X$. 

In the theory of sufficiency of jets the aim is to approach a smooth map by its family of Taylor polynomials up to some sufficient degree
depending on the kind of equivalence considered for maps (right, left, or $V$ equivalence).   
Embedding a germ into a convenient deformation is a seminal 
way of thinking for R. Thom that has been successfully achieved in his cobordism theory or in his works on regular stratifications providing 
regular trivializations.  
We could multiply examples in this spirit, old ones as well as recent ones (from recent developments in general deformation theory itself for instance
\footnote{As recalled by Kontsevich and Soibelman, Gelfand quoted that ``any area of mathematics is a kind of deformation theory", see \cite{KonSoi}.}), 
but in this introduction we will focus only on two specific examples, that will be developed thereafter: the Milnor fibration and the Lipschitz-Killing curvatures.

 \subsection*{The Milnor fibre of a complex singularity}
The first of these two examples is provided by the Milnor fibre of a complex singular analytic hypersurface  germ
$$f: (\CC^n,0)\to (\CC,0).$$ We will assume for simplicity that this singularity is an isolated one, that is to say that we will assume that $0=f(0)$ is the only critical value of $f$, at least locally around $0$. We denote $B_{(0,\eta)}$ the open ball of radius $\eta$, centred at $0$ of the ambient space depending upon the context. 
Now, for $ \eta>0$ small enough and $0<\varrho \ll \eta$, the family $(f^{-1}(\eps)\cap 
B_{(0,\eta)})_{0<\vert \eps \vert <\varrho}$ is a smooth bundle, with projection $f$, over the punctured disc $B_{(0,\varrho)}\setminus\{ 0\}\subset \CC$.
The topological type of a fibre 
$$X_\varepsilon:=f^{-1}(\eps)\cap B_{(0,\eta)}$$
 does not depend on the choice of $\eps$, and the homotopy type of this fibre
is the homotopy of a finite CW complex of dimension $n-1$, the one of a bouquet of $\mu$ spheres $S^{n-1}$, where $\mu $ is called the Milnor number of the fibration (see \cite{Mil}). On the other hand, the special singular fibre 
$$X_0:=f^{-1}(0)\cap B_{(0,\eta)}$$
 is contractible, as a germ of a semialgebraic set. It follows that the family $(X_\varepsilon)_{0<\vert \eps \vert < \eta }$ approximates
the singular fibre $X_0$ up to $\mu$ cycles that vanish as $\eps$ goes to $0$. The number $mu$ of these cycles appears as an analytic invariant of the germ of the hypersurface $f$ that is geometrically embodied on the nearby fibres $X_\varepsilon $ of the deformation on the singular fibre (see also \cite{De}). 
In \cite{Tei73}, B. Teissier embedded the Milnor number $\mu$ in a finite sequence of integers  in the following way. For a generic vector space $V$ of $\CC^n$ of dimension $n-i$, the Milnor number of the restriction of $f$ to $V$ does not depend on $V$ and is denoted $\mu^{(n-i)}$. In particular $\mu=\mu^{(n)}$ and therefore the sequence $\mu^{(*)}:=(\mu_0, \cdots, \mu^{(n-1)}, \mu^{(n)})$ gives a multidimensional version of $\mu$. 

We can consider other invariants attached to the Milnor fibre of $f$, also extending the simple invariant $\mu$: the Lefschetz numbers of the iterates of the monodromy of the Milnor fibration, that we introduce now in order to fix notations in the sequel. 
The Milnor fibre $X_\varepsilon$ may be endowed with an isomorphism $M$, the monodromy of the Milnor fibre, defined up to homotopy and that induces in an unambiguous way an automorphism, also denoted $M$, 
on the cohomology group $H^\ell(X_\varepsilon,\CC)$ 
$$M: H^\ell(X_\varepsilon,\CC)\to  H^\ell(X_\varepsilon,\CC), \ell=0,\cdots, n-1. $$
For the $m$-th iterate $M^m$ of $M$, for any $m\ge 0$, one finally defines  the Lefschetz number $\Lambda(M^m)$ of $M^m$  by 
$$ \Lambda(M^m):= \sum_{i=0}^{n-1}(-1)^itr(M^m, H^i(X_\varepsilon,\CC)),$$
where $tr$ stands for the trace of endomorphims. 
Note that $ \Lambda(M^0)=\chi(X_0)=1+(-1)^{n-1}\mu$
and that the eigenvalues of $M$ are roots of unity (see for instance \cite{SGA7}).

 A more convenient deformation of $f^{-1}(0)$ than the family  $(f^{-1}(\eps)\cap 
B_{(0,\eta)})_\varepsilon$, at least for the practical computation of the topological invariants we just have introduced, is provided by an adapted resolution of the singularity of $f$ at $0$.  To define such a resolution and fix the notations used in Section \ref{Section Mot},  let us consider 
$\sigma : (M,\sigma^{-1}(0))\to 
(\CC^n,0)$  a proper birational map which is an isomorphism over the (germ of the) complement
 of $f^{-1}(0)$ in $(\CC^n,0)$, such that $f\circ \sigma$ and  the jacobian determinant
$\jac\ \sigma$ are normal crossings and $\sigma^{-1}(0)$ is a union of components of the divisor $(f\circ \sigma)^{-1}(0)$. We denote by 
$ E_j$, for $ j\in \cal J$, the irreducible components of 
$(f\circ \sigma)^{-1}(0)$ and
assume that $E_k$ are the irreducible components of 
$\sigma^{-1}(0)$ for $k\in \cal K \subset \cal J$. For $j\in \cal J$ 
we denote by $N_j$ the multiplicity 
 $mult_{E_j}f\circ \sigma $ of 
$f\circ \sigma$ along $E_j$ and for $k\in \cal K$ by $\nu_k$ the number $\nu_k=
1+mult_{E_k}\jac \ \sigma$. For any $I\subset \cal J$, we set
$E^0_I=(\bigcap_{i\in I} E_i)\setminus (\bigcup_{j\in \cal J\setminus I}E_j)$.
The collection $(E_I^0)_{I\subset \cal J}$ gives a canonical stratification of the divisor
$f\circ \sigma=0$, compatible with $\sigma=0$ such that in some affine open subvariety $U$ in $M$ we have $f\circ \sigma (x)=u(x) \prod_{i\in I} x_i^{N_i}$, where
$u$ is a unit, that is to say a rational function which does not vanish on $U$, and $x=(x',(x_i)_{i\in I})$ are local coordinates.
Now the nearby fibres $X_\varepsilon$ are isomorphic to their lifting $\widetilde X_\eps:=\sigma^{-1}(f^{-1}(\eps)\cap B_{(0,\eta)})$ in $M$ and the family $(\widetilde X_\eps)_{0<\vert \eps\vert < \eta}$
approximates the divisor $\widetilde X_0:=\sigma^{-1}(f^{-1}(0)\cap B_{(0,\eta)})$. Of course the geometry of $\widetilde X_0$ has apparently nothing to do with the geometry of our starting germ $(f^{-1}(0),0)$, but as the topological information concerning the singularity of $f: (\CC^n,0)\to (\CC,0)$ is carried by the nearby fibres family $(X_\varepsilon)_{0<\vert \eps\vert < \eta}$, all this information is still encoded in the family $(\widetilde X_\eps)_{0<\vert \eps\vert < \eta}$, and the discrete data $N_j, \nu_k$ although depending on the choice of the resolution, may be combined in order to explicitly compute invariants of the singularity. Not only $\mu$, the most elementary of our invariants, may be computed
in the resolution, but also more elaborated ones such that the Lefschetz numbers of the iterates of the monodromy of the singularity. 
Indeed, by \cite{ACA} we have the celebrated A'Campo formulas
$$ \Lambda(M^m)=\sum_{ i\in \cal K, \ N_i/m } N_i \cdot \chi(E_{\{i\}}^0),  \ m\ge 0$$
and in particular
$$ 1+(-1)^{n+1}\mu=\chi(X_0)=\Lambda(M^0)=\sum_{ i\in \cal K} N_i \cdot \chi(E_{\{i\}}^0). $$

\begin{rem}\label{MilnorClosure}
Denoting $\bar X_0$ the closure of the Milnor fibre of $f:(\CC,0)\to (\CC,0)$, 
since the boundary $\bar X_0\setminus X_0$ is a compact smooth manifold with odd dimension, we have $\chi(\bar X_0\setminus X_0)=0$, and in 
particular $\chi(X_0)=\chi(\bar X_0)$. This is why, in the complex case and for topological considerations at the level of the Euler-Poincar\'e characteristic, the issue of the open  
or closed nature of  balls is not so relevant. In contrast, in the real case, this issue really matters. 
\end{rem}

\subsection*{Metric invariants coming from convex geometry}
The second main example of invariants arising from a deformation that we aim to emphasize and develop here, comes from convex geometry. In this case, starting from  a compact convex set of $\RR^n$, it is usual to approximate this set by its family
of $\eps$-tubular neighbourhoods, $\eps>0$, since those neighbourhoods remain convex and generally have a more regular shape than the original set. This method is notably used in 
\cite{St} to generate a finite sequence of metric invariants attached to a compact convex polytope $P$ 
(the convex hull of a finite number of points) in $\RR^n$ (actually in 
$\RR^2$ or $\RR^3$ in \cite{St}). It is established in \cite{St}
that the volume of the tubular neighbourhood of radius $\eps\ge 0$ of $P$, 
$$ T_{P,\eps}:=\bigcup_{x\in P} \bar{B}_{(x,\eps)},$$
where   $\bar{B}_{(x,\varepsilon)}$ is the closed ball of $\RR^n$  centred at $x$ with radius $\varepsilon$,
is a polynomial in $\eps$ with coefficients $\Lambda_0(P), \cdots, \Lambda_n(P)$ depending only on $P$ and being invariant 
under  isometries of $\RR^n$. We have  

$$ \forall \eps\ge 0,   \ Vol_n( T_{P,\eps})=
\sum_{i=0}^n \alpha_i\Lambda_{n-i}(P)\cdot \eps^i,
\eqno(1)$$
It is convenient to normalize the coefficients $\Lambda_i(P)$ by the introduction,  
in the equality $(1)$ defining them,  of the $i$-volume $\alpha_i$ of the $i$-dimensional unit ball.

When $\eps=0$ in this formula, one gets $\Lambda_n(P)=Vol_n(P)$.
On the other hand, denoting $\delta=\max\{\vert x-y\vert; x,y\in P\} $ the diameter of $P$, for any $x\in P$, the inclusions $B_{(x,\eps)}\subset  T_{P,\eps} \subset 
B_{(x,\eps+\delta)}$ show that $Vol_n( T_{P,\eps})
\build{\sim}_{\eps \to \infty}^{}
\alpha_n\cdot \eps^n$ and thus $\Lambda_0(P)=1$.
Denoting the Euler-Poincar\'e characteristic by $\chi$ and having in mind further generalizations, the relation $\Lambda_0(P)=1$ has rather to be considered as $\Lambda_0(P)=\chi(P)$.
A direct proof of $(1)$ leads to an expression of the other coefficients $\Lambda_i(P)$ in terms of some geometrical data of $P$.
To give this proof, we set now some notation. 

For $P$ a polytope in $\RR^n$ of dimension $n$, generated by 
$n+1$ independent points, an affine hyperplane in $\RR^n$ 
generated by $n$ of these points is called  a facet of $P$. The normal vector to a facet $F$ of $P$ is the unit vector orthogonal to $F$ and pointing in the half-space defined by $F$ not containing $P$. 
For $i\in \{0, \cdots , n-1 \}$, a { $i$-face } of $P$ is 
the intersection of $P$ with $n-i$ distinct facets of $P$.
We denote $\cal F_{i}(P) $ the set of $i$-faces of $P$. 
By convention $\cal F_n(P)=\{  P \}$.
For $x\in P$ one consider $F_x$, the unique face of $P$ of 
minimal dimension containing $x$.
If $x\in \partial P$ (the boundary of $P$), 
the {  normal exterior cone} of $P$ at $x$, denoted  $C(x,P)$, is the $\RR_+$-cone of $\RR^n$ generated  by the normal vectors to the facets of $P$ containing $x$. By convention
$C(x,P)=\{ 0 \}$, for $x\in P \setminus \partial P$.

We note that 
for $F_x$ of dimension $i\in \{0, \cdots, n-1 \}$,
$C(x,P)$ is a $\RR_+^\times$-invariant cone of $\RR^n$ of  dimension $n-i$. Furthermore, for any
$y \in F_x$, $C(x,P)=C(y,P)$. One thus defines 
$C(F,P)$, the {exterior normal cone of $P$ along a face $F$ of $P$}, by
$C(x,P)$, where $x$ is any point in $F$. One has
$C(P,P)=\{ 0 \}$.

For $P$ a degenerated polytope of $\RR^n$, that is to say that the affine subspace $[P]$ of $\RR^n$ generated by $P$ is of dimension
$<n$, one denotes $C_{[P]}(x,P)$ the exterior normal cone of $P$ at $x$ in $[P]$, since $P$ is of maximal dimension in $[P]$.
With this notation,
the exterior normal cone of $P$ at $x$ in $\RR^n$,
denoted $C_{\RR^n}(x,P)$, or simply $C(x,P)$ when no confusion is possible, is defined by $C_{[P]}(x,P)\times [P]^\perp $. We finally define $C(F,P)$, for $P$ general, as $C(x,P)$ for any $x\in F$.
The exterior normal cone of $P$ depends on the ambient space in which we embed $P$, but we now define an intrinsic measure attached to the normal exterior cone, the exterior angle.

\begin{defi}
Let $P$  be a  polytope of  $\RR^n$ and $F\in   \cal F_i(P)$. One defines the exterior angle $\gamma(F,P)$ of $P$ along $F$  (see fig.1),  by
$$\gamma(F,P):=  {1\over \alpha_{n-i}}\cdot 
Vol_{n-i}(C(F,P)\cap \bar{B}_{(0,1)})=
Vol_{n-i-1}(C(F,P)\cap S_{(0,1)}).  $$
By convention  $\gamma(P,P)=1$.
 \end{defi}

\vskip0cm
\vglue0cm
 \vglue0cm
\vskip5,5cm
\hskip1cm\includegraphics{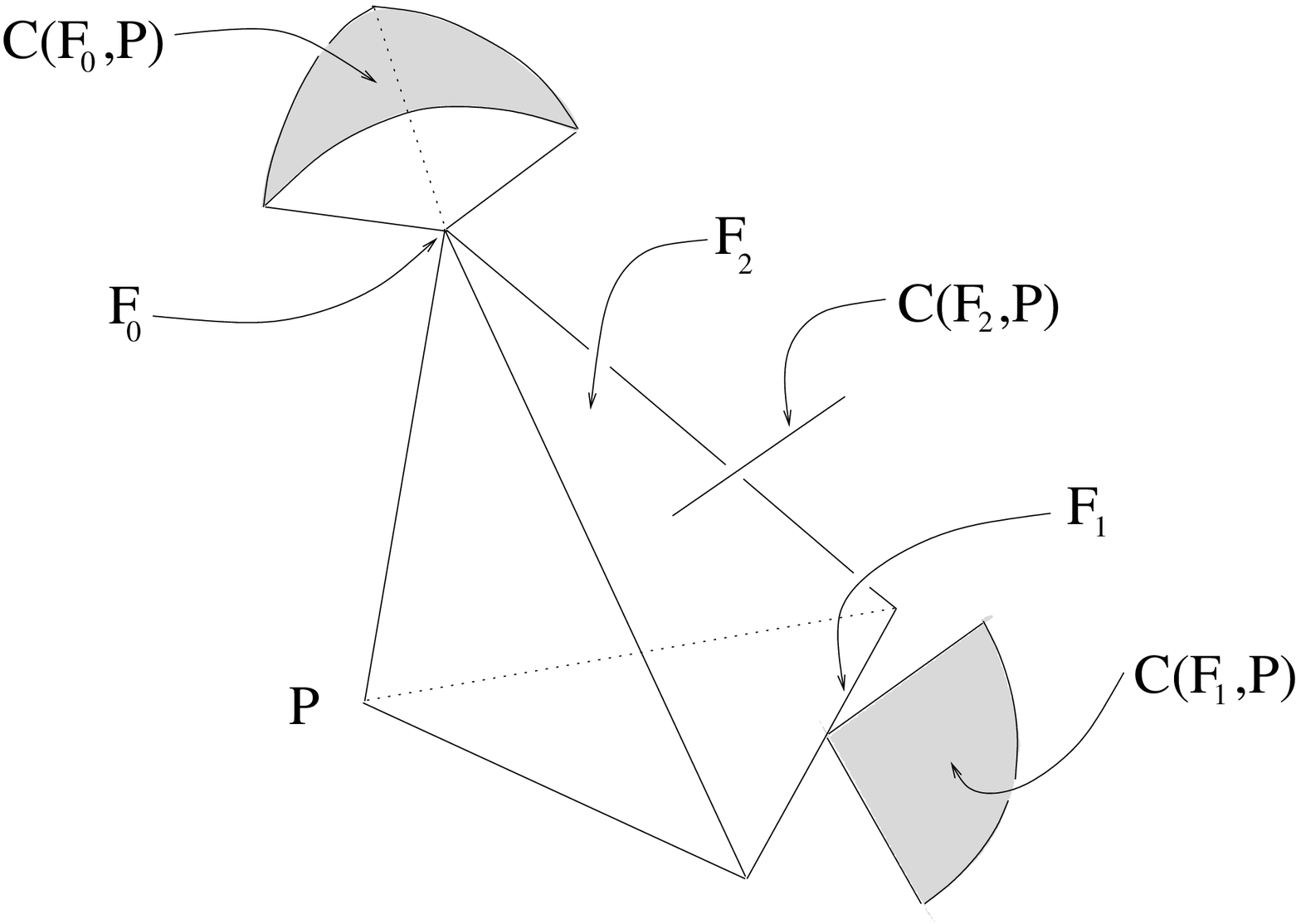}
\vglue0cm
 \vskip0,3cm
 \centerline{\bf fig.1}
 \vskip3mm
With the definition of the exterior angle, the proof of $(1)$ is trivial.
 
\begin{proof}[Proof of equality $(1)$] 
We observe that 
 $$Vol_n(   T_{P,\eps})= \sum_{i=0}^n \alpha_i\cdot \eps^{n-i}
\sum_{F\in {\cal F}_i(P)} Vol_i(F)\cdot \gamma(F,P).$$
In particular  
$$\Lambda_i(P)=\sum_{F\in {\cal F}_i(P)} 
Vol_i(F)\cdot \gamma(F,P). \eqno(2) $$ 

\end{proof}
The equality $(2)$ shows how the invariant $\Lambda_i(P)$ 
captures the concentration $\gamma(F,P)$ of the curvature of the family $( T_{P,\eps})_{\eps>0}$ along 
the $i$-dimensional faces of $P$ as $\eps$ goes to $0$.

In the general convex case and not only in the convex polyedral case, the equality $(1)$ still holds, defining invariants $\Lambda_i$ on the set $\cal  K^n$ of convex sets of $\RR^n$. A proof of this equality by approximation of a convex set by a sequence of polytopes is given in   \cite{Sch3}, section 4.2. Another proof is indicated in 
  \cite{Fed1} (3.2.35) and \cite{LanShi}, using the Cauchy-Crofton  
formula, a classical formula in integral geometry, that we recall here.  
 
\begin{CCF}[\cite{Fed1} 5.11, \cite{Fed3} 2.10.15, 3.2.26, \cite{San} 14.69]\label{CCF}
 Let  $A\subset \RR^n$ a $(\cal H^d,d)$-rectifiable set, where $\cal H^d$ is the $d$-dimensional Hausdorff measure. We have 

$$Vol_d(A)={1 \over \beta(d,n)} \int_{\bar{P} \in \overline {G}(n-d,n)} 
Card(A\cap \bar P)\  d\overline{\gamma}_{n-d,n}(\bar{P}), \eqno(\cal C   \cal C)$$
with $\overline {G}(n-d,n)$   the Grassmannian of 
$(n-d)$-dimensional  affine
planes $\bar P$ of $\RR^n$,  $\overline{\gamma}_{n-d,n}$ its canonical measure and 
denoting $ \Gamma$  
the Euler function, $\beta(d,n)$   the universal constant
 $\Gamma({n-d+1\over 2})\Gamma({d+1\over 2}) /
\Gamma({n+1\over 2}) \Gamma({1\over 2})  $. 
\end{CCF}

One can now prove equality $(1)$ in the general compact convex case. 

 \begin{proof}[Proof of equality $(1)$ in the convex case]
 We proceed by induction on the dimension of the ambient space is which our convex compact set $K$ lies. 
If this dimension is $1$, 
formula $(1)$ is trivial, and if this dimension is $n>1$, 
one has  
$$ Vol_n(  T_{K,\eps})= Vol_n(K)+\int_{r=0}^\eps Vol_{n-1}(K^r) \ dr,$$
where $K^r$ is the set of points in $ \RR^n$ at distance $r$ of $K$.
We compute $Vol_{n-1}(K^r)$ using 
the Cauchy-Crofton formula.

Noting that  $Card(\overline L \cap K^r)=2 $ or $Card(\overline L \cap K^r)=0$, 
up to a $ \overline \gamma_{1,n}$-null subset of $\overline{G}(1,n)$, we obtain by definition of 
$\bar \gamma_{1,n}$ 
$$ \hskip0mm Vol_n(  T_{K,\eps})=$$
$$ Vol_n(K)+\int_{r=0}^\eps
{2\over \beta(1,n)}\int_{H\in G(n-1,n)} Vol_{n-1}(\pi_H(K^r)) 
\ d \gamma_{n-1,n}(H)  \ dr,$$
where $G(n-1,n)$ is the Grassmannian of $(n-1)$-dimensional vector subspace of
$ \RR^n$ equipped with its canonical measure $\gamma_{n-1,n}$ invariant under the action of
$O_n( \RR)$ and $\pi_H$ is the orthogonal projection onto $H\in G(n-1,n)$.

\vskip0cm
\vglue0cm
 \vskip6,0cm
\hskip0cm\includegraphics{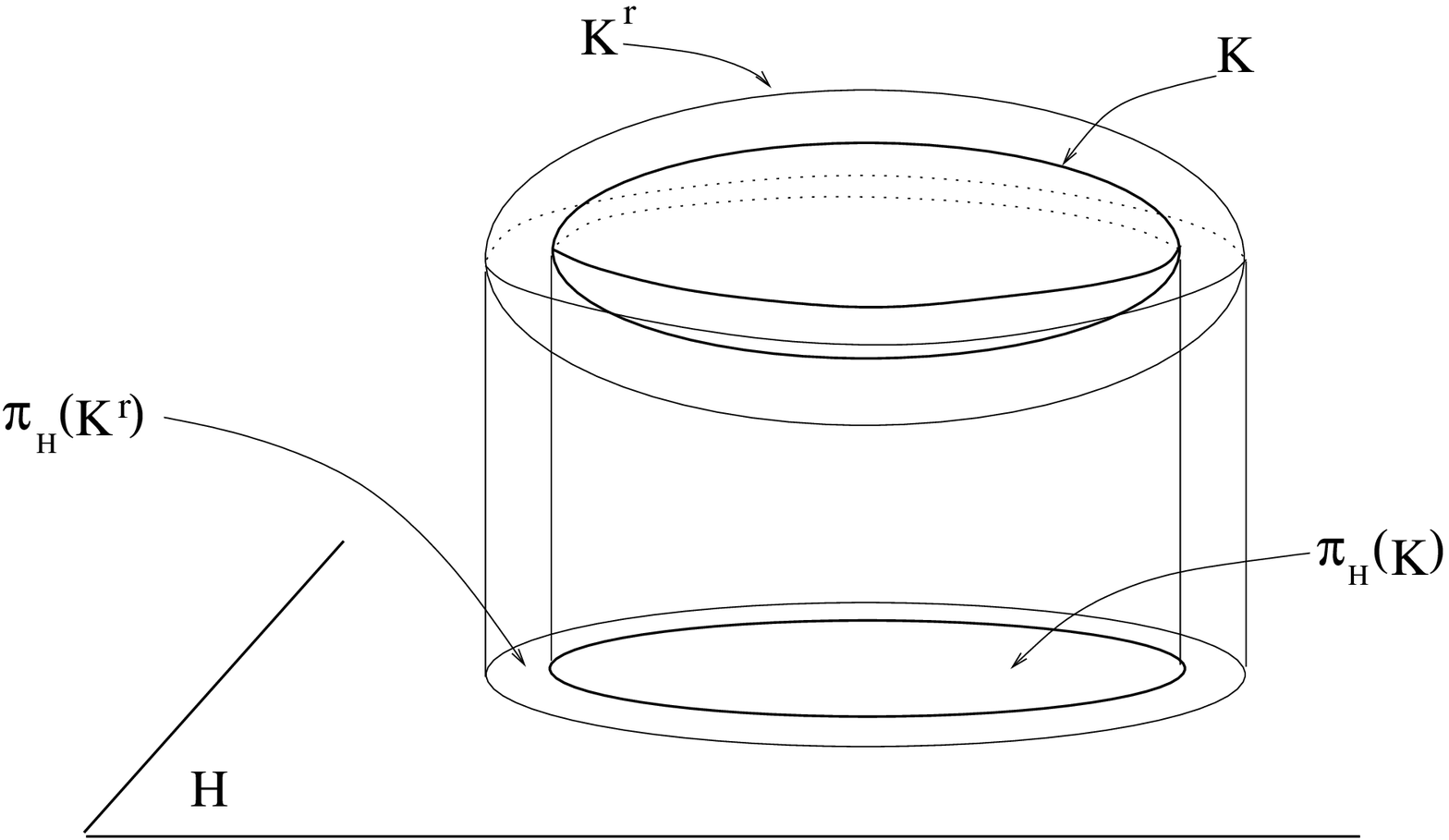}
\vglue0cm
 \vskip0,5cm
 \centerline{\bf fig.2}
 \vskip0,5cm

 By induction hypothesis, the expression of the volume of the tubular neighbourhood of radius $r$ of the convex sets of $ \RR^{n-1}$ is a polynomial in $r$. Since $\pi_H(K^r)$ is
$  T_{\pi_H(K),r}$ in  $H$, we have
$$ \hskip0mm Vol_n(  T_{K,\eps})= Vol_n(K)$$
$$+
{2\over \beta(1,n)}\int_{r=0}^\eps \int_{H\in G(n-1,n)} 
\sum_{i=0}^{n-1} \alpha_i\Lambda_{n-1-i}(\pi_H(K))\cdot r^i   
\ d \gamma_{n-1,n}(H)  \ dr$$
$$=Vol_n(K)+ {2\over \beta(1,n)}  
\sum_{i=0}^{n-1} {\alpha_i\over i+1}\cdot  \eps^{i+1}  
\int_{H\in G(n-1,n)} \Lambda_{n-1-i}(\pi_H(K)) 
\ d \gamma_{n-1,n}(H).  $$
 \end{proof}

In \cite{St}, the formula $(1)$ is also proved for $C^{2+}$ surfaces, giving a hint for a possible extension of this formula to the smooth case. This extension is due to H. Weyl, who proved in 
\cite{Wey} the following statement  (see also \cite{LanShi}). 
 
\begin{theorem}[Weyl's tubes formula]
Let  $X$ be a smooth compact submanifold of   $ \RR^n$ of dimension $d$. Let  $\eta_X >0$ such that for any $\eps$, $0<\eps \le \eta_X$, 
for any $y\in   T_{X,\eps}$,
there exists a unique $x\in X$ such that $y\in x+(T_xX)^\perp$. 
Then for any $\eps\le \eta_X$
$$Vol_n(  T_{X,\eps})= \sum_{i=0}^{[d/2]} \alpha_{n-d+2i}
\Lambda_{d-2i}(X)\cdot 
\eps^{n-d+2i},$$
where the $\Lambda_k(X)'s$ are invariant under isometric embeddings of $X$ into Riemannian manifolds. \end{theorem}  

When, on the other hand, $X$ is a non convex union of two polytopes 
$P,Q$,  $Vol_n(  T_{X,\eps})$ is no more necessarily a polynomial in $\eps$. 
For example for $X_1=P\cup Q$ where 
$P=\{(0,0)\}\subset  \RR^2$ and
$Q=\{(0,2)\} \subset  \RR^2$, and for $1 \le\eps \le 2$.
In the same way, when $X$ is a singular set, for any $\eps >0$, $Vol_n(  T_{X,\eps})$ is not necessarily a polynomial in $\eps$. 
For example for 
$X_2=\{(x,y)\in  \RR^2; x\ge 0,  \ \ (x^2+y^2-1)(x^2+(y-2)^2-1)=0\}$, 
for any sufficiently small $\eps>0$, 
$Vol_2(  T_{X,\eps})=\displaystyle (1+\eps)^2\arccos({1\over 1+\eps})
-\sqrt{\eps^2+2\eps}$ (see fig.3).
 
\begin{rem}
Nevertheless by \cite{CoLiRo}  we know  that for $X$ a subanalytic 
subset of $ \RR^n$,
$Vol(  T_{X,\eps})$ is a polynomial in subanalytic functions with variable  $\eps$ and the logarithms of these functions, and thus that it is defined in some
$o$-minimal structure over the reals. 
\end{rem}

\begin{rem}
The grey areas $\Sigma_1$ et $\Sigma_2$ in figure 3, 
counted with multiplicities $1$
in $Vol_2(  T_{X_i,\eps})$ have non polynomial contributions. 
But when these areas are counted with multiplicity 2, on one hand, with this modified computation for 
 $Vol_2(  T_{X_1,\eps})$, we obtain the sum of the areas 
 of two discs of radius $\eps$ centred at  $P$ and $Q$ and, on the other hand, with this modified computation for $Vol_2(  T_{X_2,\eps})$, we obtain twice the volume of the tubular neighbourhood of radius
$\eps$ of half a circle minus the volume of a ball of radius $\eps$.
\end{rem}

\vskip0cm
 \vglue0cm
\vskip7,0cm
\hskip2,4cm\includegraphics{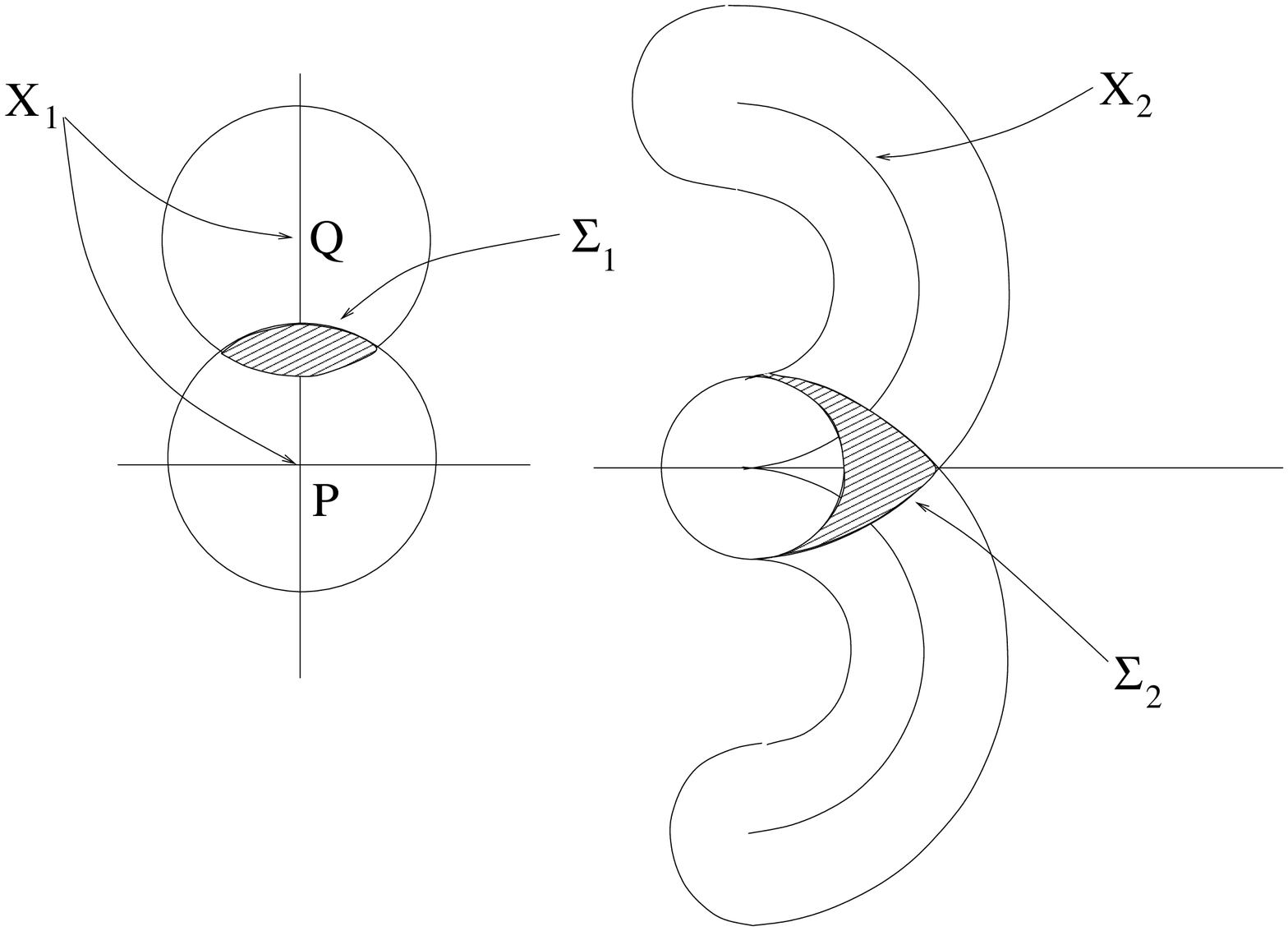}
\vglue0cm
 \vskip0,5cm
 \centerline{\bf fig.3}
  \vskip0,5cm

In conclusion, a multiple contribution of the volume of the grey areas provides two polynomials in $\eps$.
Moreover, we observe that for $j=1,2$:
\vskip2mm
- $\forall x\in \Sigma_j$: 
$2=\chi(X_j\cap \bar{B}_{(x,\eps)})$, 
\vskip1mm
- $\forall x\in   T_{X_j,\eps}\setminus \Sigma_j$:
$1=\chi(X_j\cap \bar{B}_{(x,\eps)}),$ 
\vskip1mm
- $\forall x\in  \RR^2\setminus   T_{X_j,\eps}$ : $\chi(X_j\cap
\bar{B}_{(x,\eps)})=0$. 
\vskip2mm
 
It follows that for j=1,2
$ \displaystyle \int_{x\in \R^2} \chi(X_j\cap \bar{B}_{(x,\eps)})\ dx=
\int_{x\in T_{X_j,\eps}} \chi(X_j\cap \bar{B}_{(x,\eps)})\ dx$
is a polynomial in $\eps$.
 These examples are in the scope of a general fact: the formula $(1)$
can be generalized to compact sets definable in some o-minimal structure over the reals by the formula $(1')$ given below.
 
\begin{theorem}[\cite{Fu1}, \cite{Fu2}, \cite{Fu3}, \cite{Fu4}, \cite{Fu5}, \cite{Fu6}, \cite{Fu7}, \cite{BeBr1}, \cite{BeBr2}, 
\cite{BrKu}]\label{LKpoly}
Let  $X$ be a compact subset of $\RR^n$ definable in some o-minimal structure over the ordered real field.  
There exist constants $\Lambda_0(X), \cdots, \Lambda_n(X)$ such that for any  $\eps\ge 0$ 
$$ \int_{x\in {T_{\eps,X}}}\chi(X\cap \bar{B}_{(x,\eps)})\ dx=
\sum_{i=0}^n \alpha_i\Lambda_{n-i}(X)\cdot \eps^i.
\eqno(1')$$
The real numbers  $\Lambda_i(X), i=0,\cdots, n$ are called the
Lipschitz-Killing curvatures of $X$, they only depend on definable isometric embeddings of $X$ into euclidean spaces. Moreover, we have 
$$ \Lambda_i(X)= \int_{\bar P \in \bar G(n-i,n)} \chi(X \cap \bar P) \ \
{d \bar\gamma_{n-i,n}(\bar P) \over \beta(i,n)}. \eqno(2')$$ 
\end{theorem}

\begin{remabout}

In Theorem \ref{LKpoly} we assume the set $X$ compact, although for $X$ bounded but not compact the equality $(1')$ together with $(2')$ 
is still true with $\chi$ the Euler-Poincar\'e characteristic with compact support,  usually considered for non-compact definable sets. 
This characteristic is additive and multiplicative and defined by any finite cell decomposition $\cup_i C_i$ of $X$ by 
$\chi(X)=\sum_i (-1)^{\dim(C_i)}$ (see \cite{Dri}, p. 69). For simplicity, in what follows we will still consider the compact case. 

\begin{rem}
The formula $(1')$ is clearly a generalization to the non convex case of the formula $(1)$, since for $X$ compact convex, for any  $\eps\ge 0$, for any  
$x\in   T_{X,\eps}$, $\chi(X\cap \bar{B}_{(x,\eps)})=1$, thus 
 $ \int_{x\in { T_{\eps,X}}}\chi(X\cap \bar{B}_{(x,\eps)})\ dx=
Vol_n(  T_{X,\eps})$. In the same way $(1')$ generalizes Weyl's tube formula to the singular case, since for $X$ 
smooth, there exists $\eta_X>0$ such that for any  $\eps, 0<\eps<\eta_X$, for any  $x\in   T_{X,\eps}$, $\chi(X\cap \bar B_{(x,\eps)})=1$ and 
 $ \int_{x\in { T_{\eps,X}}}\chi(X\cap \bar B_{(x,\eps)})\ dx=
Vol_n(  T_{X,\eps})$.
\end{rem}

The formula $(1')$ comes from a more general cinematic  formula
(see \cite{BrKu}, \cite{Fu5}). For $X$ and $Y$ two  definable sets of $ \RR^n$  
$$ \int_{g\in G}\Lambda_k(X\cap g\cdot Y)\ dg=
\sum_{i+j=k+n} c_{n,i,j}\cdot\Lambda_{i}(X)\cdot \Lambda_{j}(Y)
$$
with $G$ the group of isometries of $ \RR^n$ and $c_{n,i,j}$ universal constants. 

 The expression of $\Lambda_i$ given by $(2')$ generalizes to the definable case the representation formula $(2)$ of $\Lambda_i$ given in the polyhedral case. Furthermore, from $(2')$ we get the following characterization of  $\Lambda_0$ and $\Lambda_d$,  $d=\dim(X)$, already obtained from $(1)$ in the compact convex case $$\Lambda_0(X)=\chi(X)$$ and, using the Cauchy-Crofton formula, $$\Lambda_d(X)=Vol_d(X).$$ Finally, 
$$\Lambda_{d+1}(X)=\cdots=\Lambda_n(X)=0,$$ for $d<n$.

The last remark made now here is a remark that, having in mind geometric measure theory, we are eager to address: the Euler-Poincar\'e  
characteristic being additive for definable sets (see \cite{Dri})
the equality  $(1')$ or $(2')$ shows that the $ \Lambda_i$'s are additive invariants of   definable sets, in the following sense 
 
$$\forall i\in \{0, \cdots, n\},  \ \Lambda_i(X\cup Y)
= \Lambda_i(X)+ \Lambda_i(Y)- \Lambda_i(X\cap Y), $$
for any   definable sets  $X$ and $Y$ of $ \RR^n$. 
\end{remabout}

We have now in hand two kinds of deformation of a singular set. When this set is an analytic isolated hypersurface singularity, we may consider its Milnor fibration, providing in particular as invariant the Milnor number of the singularity, and in the more general definable case, we have recalled in details the notion of Lipschitz-Killing curvatures, coming from the deformation family provided by the tubular neighbourhoods.   
It is worth noting that in these two cases, the deformations considered lead to additive invariants attached to the given germ. 
 
In what follows, we explain how to localize the Lipschitz-Killing curvatures in order to attach to a singular germ a finite sequence of additive invariants (Section 1). 
We will then explain how all our local invariants are related and how such kind of relation illustrates,  in the very general context of definable sets over the reals, the emergence of a well-known principle in convex geometry wherein additive invariants (with some additional properties) may not be so numerous (Section 2). 
Finally, we will stress the fact that the additive nature of an invariant coming from a deformation allows us to compute this invariant in some adapted  scissors ring via some generating zeta function capturing the nature of the deformation. This is the point of view underlies the work of  Denef and Loeser (Section 3).

\begin{nota}
As well as in this introduction, in the sequel, $B_{(x,r)}$, $\bar{B}_{(x,r)}$ and $S_{(x,r)}$ are respectively the open ball, the closed ball and
the sphere centred at $x$ and with radius $r$ of the real or the complex vector spaces $\RR^n$ or 
$\CC^n$. 
If necessary, to avoid confusion,  we emphasize the dimension $d$ of the
ambient space   to which the ball belongs by denoting $B^d_{(x,r)}$. 
 Definable means definable in some given $o$-minimal structure expanding the ordered real field $(\RR,+,-,\cdot,0,1,<)$ 
(see \cite{Cha}, \cite{Dri}).
\end{nota}

\tableofcontents

\section{Local  invariants from the tubular neighbourhoods deformation}
  
The invariants $\Lambda_0, \cdots, \Lambda_n$ defined in the introduction for compact definable sets (or at least bounded definable sets) may be extended to non bounded definable sets as well as they may be localized in order to be attached to any definable germ $(X,0)$. The extension of the invariants $\Lambda_i$
to non bounded definable sets has been proposed in \cite{Du08}. 
These two possible extensions are similar; they essentially use the fact that near a given point or near infinity the topological types of affine sections of a definable set are finite in number. As we are mainly interested in local singularities, we explain in this section how to localize the sequence  $(\Lambda_0, \cdots, \Lambda_n)$ at a given point.

For this goal, let us consider $X\subset \RR^n$ a compact definable set. We assume that $0\in X$ and we denote by $X_0$ the germ of $X$ at $0$, $d$ its dimension. 
Representing elements $\bar P$ of the Grassmann manifold $\bar G(n-i,n)$ of $(n-i)$-dimensional affine subspaces of $\RR^n$
 by pairs $(x, P) \in \RR^n \times G(i,n)$, where
$x \in P$, and $\bar P$ is the affine subspace of $\RR^n$ orthogonal to $P$ at $x$,
the measure $\bar \gamma_{n-i,n}$ on $\bar G(n-i,n)$ is the image 
through this representation of the product $m \otimes \gamma_{i,n}$, where $m$ is the Lebesgue measure on $P$ and $P$ is identified with $\RR^i$.
It follows by formula (2') that $\Lambda_i$ is $i$-homogeneous, that is to say 
$\Lambda_i(\lambda\cdot X)=\lambda^i\Lambda_i(X)$, 
for any $\lambda\in   \RR_+^*$. In consequence, it is natural to consider the asymptotic behaviour of  
$ \displaystyle \frac{1}{\alpha_i}\Lambda_i(\frac{1}{\varrho}\cdot (X \cap \bar{B}_{(0,\varrho)}))=\frac{1}{\alpha_i\varrho^i}\Lambda_i
(X\cap \bar{B}_{(0,\varrho)})$,
as $\varrho\to 0$, in order to obtain invariants attached to the germ   $X_0$ of $X$ at $0$. 
Using standard arguments for the definable family 
$$({1\over \varrho}\cdot (X\cap \bar B_{(0,\varrho)})\cap \bar P))_{(\varrho, \bar P)\in \RR_+^*\times \bar G(n-i,n)}$$
such as Thom-Mather's isotopy lemma or cell decomposition theorem, one knows that, for any fixed $\bar P\in  \bar G(n-i,n)$, the topological type of the family 
$({1\over \varrho}\cdot (X\cap \bar B_{(0,\varrho)})\cap \bar P))_{\varrho )\in \RR_+^*}$ is constant for $\varrho$ small enough and therefore the limit of $\chi({1\over \varrho}\cdot (X\cap \bar B_{(0,\varrho)})\cap \bar P)$ for  $\varrho \to 0$ does exist. Furthermore, 
 still by finiteness arguments proper to definable sets, the family  $(\chi({1\over \varrho}\cdot (X\cap \bar B_{(0,\varrho)})\cap \bar P))_{\bar P\in \bar G(n-i,n)}$ is bounded with respect
to $\bar P$. 

The next definition follows from these observations (see \cite{CoMe}, Theorem 1.3). 

\begin{defi}[Local Lipschitz-Killing invariants, see \cite{CoMe}]\label{Lloc}
Let  $X$ be a (compact) definable set of $\RR^n$, representing the germ $X_0$ at $0\in X$. 
The limit 
$$ \Lambda_i^{\ell oc}(X_0):=\lim_{\varrho \to 0}{1\over \alpha_i.\varrho^i} 
\Lambda_i(X \cap \bar B{(0,\varrho)}) \eqno(3)$$
exists and
the finite sequence of real numbers  $ (\Lambda_i^{\ell oc}(X_0))_{i \in \{0, \ldots
, n \} }$ is called the sequence of local Lipschitz-Killing invariants of the germ $X_0$.
\end{defi}

\begin{rem}
Another kind of localization of 
the invariants $\Lambda_i$ have been obtained and studied in \cite{BeBr2},
by considering the family  $ (X\cap S_{(0,\varrho)})_{\varrho>0}$ instead of the family 
$(X\cap \bar B_{(0,\varrho)})_{\varrho>0}$.
\end{rem}

\begin{rems}
For any $i\in \{0,\cdots, n\}$, just as $\Lambda_i$, $\Lambda_i^{\ell oc}$ is invariant 
under isometries of 
$\RR^n$ and defines an additive function on the set of definable germs at the origin of $\RR^n$.
Moreover  $\Lambda_i^{\ell oc}(X_0)=0, \ \hbox{ for } i>d$, since 
$\Lambda_i=0$  for definable sets of dimension $<i$ and for any definable compact germ 
$X_0$, $\Lambda_0^{\ell oc}(X_0)=1$, since $\Lambda_0=\chi$ and a definable germ is contractible.  
Finally, since by the Cauchy-Crofton formula $ (\cal C   \cal C)$ and $(2')$ we have $$\Lambda_d^{\ell oc}(X_0)=\lim_{\varrho \to 0}
{Vol_d(X\cap B_{(0,\varrho)})\over Vol_d(B^d_{(0,\varrho)})},\eqno(4)$$ 
we observe that $\Lambda_d^{\ell oc}(X_0)$ is by definition  the local density $\Theta_d(X_0)$ of $X_0$, and thus we have obtained,
by finiteness arguments leading to
 Definition \ref{Lloc}, the following  theorem of    Kurdyka et   Raby.
\end{rems}

\begin{cor} [\cite{KuRa}, \cite{KuPoRa}, \cite{Li}]\label{Kur}
 The local density of  definable sets of  $\RR^n$ exists at each point of $\RR^n$. 
\end{cor}

On  figure $4$ are represented the data taken into account 
in the computation of $ \Lambdal_i(X)$.

For $P\in G(i,n)$, we denote by $K^{P,\varrho}_\ell$ 
the domains of $P$ above which the Euler-Poincar\'e characteristic
of the fibres of  $\pi_{P\vert X\cap \bar B(0,\varrho)}$ is constant
and equals $\chi_\ell^{P,\varrho}\in \ZZ$.
The quantity  $\Lambda_i(X\cap \bar B^n_{(0,\varrho)})$ is then obtained as the mean value over the vector planes $P$ of the sum 
 $\sum_{\ell=1}^{\ell_P} 
\chi_\ell^{P,\varrho}\cdot Vol_i(K^{P,\varrho}_\ell)$. In particular we'd like to stress the fact that are considered in this sum 
the  volumes of the  domains $K^{P,\varrho}_\ell$ (in green on figure 4) defined
by the critical values of $\pi_P$ coming from  
the link $X\cap S_{(0,\varrho)}$. We draw attention to these green domains, far from the origin, in view of other local invariants, the polar invariants, that will be defined in the next section, and for which only the domains $K_\ell^{P,\varrho}$ close to the origin will be considered.  
\eject
 
\vskip1,3mm
$$\hskip-9,5cm  X\cap B^n_{(0,\varrho)}$$
\vskip-1,1cm
\vglue0cm
\vskip10cm
\hskip1cm\includegraphics{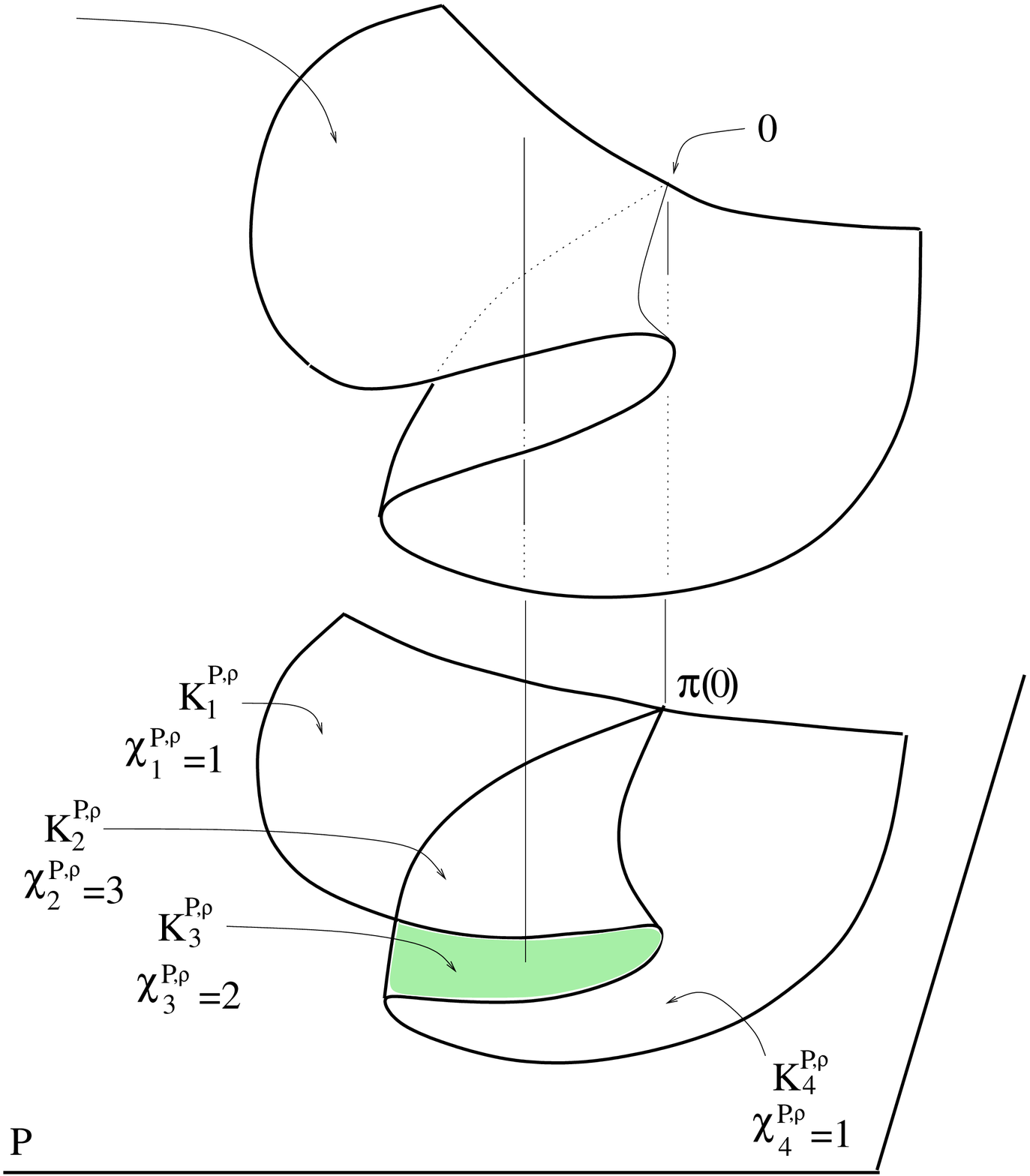}
\vglue0cm
\vskip0,5cm
\centerline{\bf fig. 4}
 \vskip0,8cm
 \noindent

\section{Additive invariants of singularities and Hadwiger principle}
In the previous section we have localized the Lipschitz-Killing invariants. We'd like now to investigate the question of how these invariants are related to classical local invariants of singularities, such as Milnor number, or Milnor numbers of generic plane sections of the singularity (for the complex case). The question of the correspondence of invariants coming from differential geometry and invariants of singularities has been tackled by several authors. In the first section \ref{section 2.1}, we briefly recall some of these works. 
We then introduce in section  \ref{section 2.2} a sequence of invariants of real singularities that is the real counterpart of classical invariants of complex singularities and we finally relate the localized  Lipschitz-Killing invariants to the polar invariants in  section \ref{section 2.3}. In section \ref{section 2.4} we recall results from convex geometry and convex valuation theory giving a strong hint, called here Hadwiger principle,  of the reason why such invariants have to be linearly related. 

\subsection{Differential geometry of complex and real hypersur\-face singularities}\label{section 2.1}
The first example  we'd like to recall of such a correspondence may be found in \cite{Lan} (see also \cite{Lan2} and \cite{LanLe}). In  \cite{Lan} R. Langevin relates the concentration of the curvature of the Milnor fibre $F_\varepsilon^\eta=f^{-1}(\varepsilon)\cap B_{(0,\eta)}$, $0<\varepsilon \ll \eta \ll 1$, as $\varepsilon$ and $\eta$ go to $0$ to the Milnor numbers $\mu$ and $\mu^{(n-1)}$ of the isolated hypersurface singularity $f:(\CC^n,0)\to (\CC,0)$. To present this relation, we need some definitions. 

For a given real smooth oriented hypersurface $H$ of $\RR^n$ and for $x\in H$, we classically define the Gauss curvature $K(x)$ of $H$ at $x$ by $K(x):=\hbox{jac}( \nu) $, where $\nu: H \to \SS^{n-1}$ is the mapping giving the unit normal vector to $H$  induced by the canonical orientation of $\RR^n$ and the given orientation of $H$. The curvature $K(x)$ can be generalized in the following way to any submanifold $M$ of $\RR^n$ of dimension $d, d\in \{0,\cdots, n-1\}$ (see \cite{Fen}). Let $x$ be a point of $M$ and denote by $N(x)\simeq \SS^{n-d-1}$ the manifold of normal vectors to $M$ at $x$, and, for $\nu\in N(x)$, by $K(x,\nu)$ the Gauss curvature at $x$ of the projection $M_\nu$ of $M$ to $\hbox{T}_xM\oplus \nu$. Note that the projection $M_\nu$ defines at $x$ a smooth hypersurface of $\hbox{T}_xM\oplus \nu$ oriented by $\nu$. 
 The mean value of $K(x,\nu)$ over 
$N(x)$ define the desired generalized Gauss curvature. 

\begin{definition}[see \cite{Fen}]\label{courbugene}
With the above notation, the curvature $K(x)$ of $M$ at $x$ is defined by 
$$K(x):=\int_{\nu\in \PP N(x)}K(x,\nu)\  d\nu $$
\end{definition} 

\begin{remarks} 
The curvature $K(x,\nu)$ is $\varepsilon^{d-n}$ times the Gauss curvature of the boundary $\partial T_{M,\varepsilon} $ of the $\varepsilon$-neighbourhood of $M$ at $x+\varepsilon\nu$ (see \cite{CheLas}). 

In \cite{Lan}, following Milnor, it is observed that for $M$ a smooth complex hypersurface of $\CC^n$, $K(x)=(-1)^{n-1} \pi \vert \hbox{jac} \gamma_\CC(x)\vert^2$, where $\gamma_\CC$ is the complex Gauss map sending $x\in M$ to the normal complex line $\gamma_\CC(x)\in \PP\CC^{n-1}$ to $M$ at $x$.

 In case $M$ is a  compact submanifold of $\RR^n$, using a so-called exchange formula (\cite{Lan2} Theorem II.1, \cite{LanShi}, \cite{Du02}) relating $\int_M K(x) \ dx$ and the mean value over generic lines $L$ in $\RR^n$ of the total index of the projection of $M$ on $L$, we obtain the Gauss-Bonnet theorem 
$$ \int_M K(x) \ dx=c(n,d)\chi(M). $$
\end{remarks}

Applying Definition \ref{courbugene} to the Milnor fibre 
$F_\varepsilon^\eta$ of an isolated hypersurface singularity $f:(\CC^n,0)\to (\CC,0)$ and using again the exchange formula, we can estimate the concentration of the curvature $K(x)$ of the Milnor fibre $F_\varepsilon^\eta$ as 
$\varepsilon$ and $\eta$ go to $0$. This value is related to the invariants $\mu$ and $\mu^{(n-1)}$ of the singularity, thanks to a result of \cite{Tei75},  by the following formula

 \begin{thm}[\cite{Lan}, \cite{Lan2}]\label{Langevin}
The curvature $K$ of a complex Milnor fibre satisfies 
 $$\lim_{\eta \to 0}\lim_{\varepsilon \to 0} \ c(n) \int_{F_\varepsilon^\eta}(-1)^{n-1} K(x)\ dx =  \mu + \mu^{(n-1)},$$
 where $c(n)$ is a constant depending only on $n$. 
 \end{thm}

This formula has been generalized to the other terms of the sequence $\mu^{(*)}$ by  Loeser in \cite{Loe84} in the following way.

 \begin{thm}[\cite{Loe84}]\label{Loeser} For $k\in \{1, \cdots, n-1\}$, we have
$$ \lim_{\eta \to 0}\lim_{\varepsilon\to 0} 
\frac{(-1)^{n-k}c(n,k)}{\eta^{2k}}\int_{F_\varepsilon^\eta}
c_{n-1-k}(\Omega_{f^{-1}(\varepsilon)}) \wedge \Phi^k
  = \mu^{(n-k)}+\mu^{(n-k-1)},$$
where $c_{n-1-k}(\Omega_{f^{-1}(\varepsilon)})$ is the $(n-1-k)$-th Chern form of  $f^{-1}(\varepsilon)$,  $\Phi$ the K\"ahler form of $\CC^n$ and as usual $c(n,k)$ a constant depending only on $n$ and $k$.
\end{thm}

A real version of these two last statements has been given by N. Dutertre, in \cite{Du02} for the real version of Theorem \ref{Langevin} and in \cite{Du08} for the real version of Theorem \ref{Loeser}.
In \cite{Du02} (see Theorem 5.6), a real polynomial germ $f:(\RR^n,0)\to (\RR,0)$ having an isolated singularity at $0$ is considered and the following equalities are given for the asymptotic behaviour of the Gauss curvature on the real Milnor fibre $F_\varepsilon^\eta=f^{-1}(\varepsilon)\cap B_{(0,\eta)}$.

 \begin{thm}[\cite{Du02}, Theorem 5.6]\label{Dutertre02}
 The Gauss curvature $K$ of the real Milnor fibre $F_\varepsilon^\eta$ have the following asymptotic behaviour
 $$\lim_{\eta \to 0}\lim_{\varepsilon \to 0^+}   \int_{F_\varepsilon^\eta} K(x)\ dx =  
 \frac{Vol(S^{n-1})}{2}\hbox{deg}_0\nabla f+ \frac{1}{2} \int_{G(n-1,n)} \hbox{deg}_0\nabla(f_{\vert{P}}) \ dP$$
 
  $$\lim_{\eta \to 0}\lim_{\varepsilon \to 0^-}   \int_{F_\varepsilon^\eta} K(x)\ dx =  
- \frac{Vol(S^{n-1})}{2}\hbox{deg}_0\nabla f+ \frac{1}{2} \int_{G(n-1,n)} \hbox{deg}_0\nabla(f_{\vert{P}}) \ dP$$
\end{thm}
 
 In \cite{Du08}, the asymptotic behaviour of the symmetric 
 functions $s_0,\cdots, s_{n-1}$ of the curvature of the real Milnor fibre are studied. For a given smooth hypersurface $H$ of $\RR^n$, the $s_i$'s are defined by $$\det(Id+tD\nu(x))=\sum_{i=0}^{n-1}s_i(x)\cdot t^i=\prod_{i=1}^{n-1}(1+k_i(x)t),$$
 where the $k_i$'s are the principal curvatures of $H$, that is to say, the eigenvalues of the symmetric morphism $D\nu(x)$. The limits
 $$\lim_{\eta \to 0}\lim_{\varepsilon \to 0}  
 \frac{1}{\eta^k} \int_{F_\varepsilon^\eta} s_{n-k}(x)\ dx $$
 are then given in terms of the mean value of  $\hbox{deg}_0\nabla(f_{\vert{P}})$ for $P\in G(n-k+1,n)$
 and for $P\in G(n-k-1,n)$ (see \cite{Du08}, Theorem 7.1).
 In particular, the asymptotic behaviour of the symmetric 
 functions $s_0,\cdots, s_{n-1}$ of the curvature of the real Milnor fibre is related to the Euler-Poincar\'e characteristic of the real Milnor fibre by the following statement. 
 
 \begin{thm}[\cite{Du08}, Corollary 7.2]\label{Dutertre08}

For $n$ odd,
 $$\chi(F_\varepsilon^\eta)=\sum_{k=0}^{(n-1)/2}c(n,k) \lim_{\eta \to 0}\lim_{\varepsilon \to 0} \frac{1}{\eta^{2k}}   \int_{F_\varepsilon^\eta} s_{n-1-2k}(x)\ dx,  $$
 and for $n$ even
 $$\chi(F_\varepsilon^\eta)=\frac{1}{2}\chi(f^{-1}(0)\cap S^{n-1}_{(0,\eta)})=\sum_{k=0}^{(n-2)/2} c'(n,k) \lim_{\eta \to 0}\lim_{\varepsilon \to 0} \frac{1}{\eta^{2k}}    \int_{F_\varepsilon^\eta} s_{n-2-2k}(x)\ dx,  $$
\end{thm}

\subsection{Local invariants of definable singular germs}\label{section 2.2}
In the present survey we aim  to relate the local Lipschitz-Killing invariants
$\Lambda_i^{\ell oc}, i=0,\cdots, n$, coming from the tubular neighbourhoods deformation, to local invariants of definable singular germs of $\RR^n$. These germs have not necessarily to be of codimension $1$ in $\RR^n$, as it is the case in the complex and real statements recalled above in Section \ref{section 2.1}. Therefore we have to define local invariants of singularities attached to definable germs of $\RR^n$ of any dimension and try to relate them to the sequence $\Lambda_*^{\ell oc}$. Furthermore those invariants have to extend, to the real setting, classical invariants of complex singularities, such as the sequence $\mu^{(*)}$ in the hypersurface case or the sequence of the local multiplicity of polar varieties in the general case. For this purpose we introduce now a new sequence $\sigma_*$ of local invariants, called the sequence of polar invariants. 

Let, as before, $X \subset \RR^n$ be a closed definable set, and assume that $X$ contains the origin of $\RR^n$ and that $d$ is the dimension of $X$ at $0$. We denote by
 $\cal C(X)$ the group of definable constructible functions on $X$, that is to say the group of definable $\ZZ$-valued functions on $X$. These functions $f$ are characterized by the existence of a finite definable partition $(X_i)$ of $X$ (depending on $f$) such that $f_{\vert X_i}$ is 
 a constant integer $n_i\in \ZZ$, for any $i$.   We denote by 
 $\cal C (X_0)$ the group of germs at the origin of functions of $\cal C (X)$. 
 For $Y\subset \RR^m $ a definable set,
$f:X\to Y$ a definable  mapping, a definable  set 
  $Z\subset X$ and  $y\in Y$, we introduce the notation $f_*( {\bf 1}_Z)(y)
:=\chi(f^{-1}(y)\cap Z)$ and we then define the following  functor from the category of definable sets to the category of groups 

$$ \begin{matrix} 
 X &  \leadsto & \cal C(X)\\
 f\downarrow \ \ \ &  &   \ \   \downarrow f_* \\
 Y & \leadsto  &  \cal C(Y)\\ 
   \end{matrix} $$
 In \cite{CoMe}, Theorem $2.6$,
it is stated that, for $f=\pi_P$ the (orthogonal) projection 
onto a generic $i$-dimensional  vector subspace $P$ of $\RR^n$, this diagram leads to the following diagram for germs
$$\begin{matrix}
 X_0 &   \leadsto & \cal C(X_0) & \\
 \pi_{P_0}\downarrow \ \ \ & & \downarrow \pi_{P_0*} \\
P_0  &   \leadsto & \cal C(P_0)\\  
\end{matrix}\eqno(5)$$
where for $Z_0\subset X_0$ and $y\in P$, $\pi_{P_0*}( {\bf 1}_{Z_0})(y)$ is defined by $\chi(\pi_P^{-1}(y)\cap Z\cap 
\bar B_{(0,\varrho)})$,
$\varrho$ being sufficiently small and $0<\Vert y\Vert\ll \varrho$.
 The existence of such a diagram for germs simply amounts to prove that a generic projection of a germ defines a germ 
 \footnote{ Let us for instance denote $X$ the blowing-up of $\RR^2$ at the origin and $x\in X$ a point of the exceptional divisor of
 $X$. We then note that the projection  of the germ $(X)_x$ on $\RR^2$, along the exceptional divisor of $X$, 
 does not define  a germ of $\RR^2$. Indeed, the projection of
 $X\cap U$, for $U$ a neighbourhood of $x$ in $X$, defines a germ at the origin of $\RR^2$ that depends of $U$.}{}
  and that, for such a projection and for any $c\in \ZZ$, the germ at $0$ of the definable set 
 $\{y\in P; \chi(\pi_P^{-1}(y)\cap Z\cap \bar B_{(0,\varrho)})=c\}$ does not depend on $\varrho$. 
 
Denoting by 
$\theta_i(\varphi)$ the integral with respect to the local density $\Theta_i$ at $0\in \RR^i$ of a germ  $\varphi:P_0 \to \ZZ$ of constructible function, that is to say  
$$\theta_i(\varphi):=\sum_{j=1}^N n_j\cdot \Theta_i(K^j_0),$$ 
when $\varphi=\sum_{j=1}^N n_j\cdot {\bf 1}_{K^j_0}$, 
for some definable germs $K^j_0\subset P_0$ partitioning $P_0$, we can define the desired  polar invariants $\sigma_i(X_0)$ of $X_0$.
\begin{definition}[Polar invariants]\label{polarinv}
With the previous notation, the polar invariants of the definable germ $X_0$ are   
 $$\sigma_i(X_0):=\int_{P\in G(i,n)} 
\theta_i( \pi_{P_0*}({\bf 1}_{X_0}) )\ d \gamma_{i,n}(P), 
\ \ i=0, \cdots, n$$
\end{definition}

\begin{remarks}

Since they are defined as mean values over generic projections, the  
  $\sigma_i$'s are invariant under the action of isometries of $\RR^n$. On the other hand the  $\sigma_i$'s define additive invariants (as well as the $\Lambda_i^{\ell oc}$'s do), since they are defined through the Euler-Poincar\'e characteristic $\chi$ and the local density $\Theta_i$, two additive invariants. 

 Observe that $\sigma_i(X_0)=0$, for $i>d$, since a general $k$-dimensional  affine subspace of $\RR^n$ does not encounter a definable set of
codimension $>k$.
We also have  $\sigma_0(X_0)=\Lambda_0^{\ell oc}(X_0)=1$, again by the local conic structure of definable sets and because $X_0$ is closed. Finally,  for $i=d$, one shows that  $$\sigma_d(X_0)=\Theta_d(X_0)=\Lambda_d^{\ell oc}(X_0)$$ 
 (we recall that by the Cauchy-Crofton formula $ (\cal C   \cal C)$ and by definition $(2')$ and $(3)$ of $\Lambda_d$  and of $\Lambda_d^{\ell oc}$, we  have $\Theta_d(X_0)=\Lambda_d^{\ell oc}(X_0)$, as already observed for Corollary \ref{Kur}). Since the relation  
 $$\sigma_d(X_0)=\Theta_d(X_0)$$
  asserts that the localization $\Theta_d$ of the $d$-volume is $\sigma_d$, that is to say, by definition of $\sigma_d$, that the localization of the volume may be computed by the mean value over (generic) $d$-dimensional vector subspaces $P\subset \RR^n$ of the number of points in the fibre of the projections of the germ $X_0$ onto the germ $P_0$, this relation appears as the local version of the global Cauchy-Crofton formula $ (\cal C   \cal C)$. We state it as follows.

\end{remarks}

\begin{LCCF}[\cite{CoCRAS}, \cite{Co} 1.16, \cite{CoMe}  3.1]\label{LCCF}
 Let $X$ be a definable subset of  $\RR^n$ of dimension $d$ (containing the origin), let  $\cal G$ be a definable subset of $G(d,n)$ on which transitively acts a subgroup $G$ of $O_n(\RR)$ and
let $m$ be a $G$-invariant measure on $\cal G$, such that 
 \par
{\bf -} the tangent spaces to the tangent cone of   de $X_0$ are in 
$\cal G$, \par
{\bf -} There exists $P^0\in \cal G$ such that   $\{g\in G; g\cdot P^0=P^0\}$  transitively acts on the $d$-dimensional vector subspace $P^0$ and 
$m(\cal G)=m(\cal G \cap \cal E_X)=1$, where $ \cal E_X$ is the generic set of $G(d,n)$ for which the localization 
$(5)$ est possible. \par
Then, we have
$$ \sigma^{\cal G}_d(X_0)=\Theta_d(X_0), \eqno(\cal C\cal C^{\ell oc}) $$
where $\sigma_d^{\cal G}$ is defined as in Definition \ref{polarinv}, but relatively to $\cal G$ and $m$.
\end{LCCF}

In the case $\cal G=G(d,n)$ and 
 $G=O_n(\RR)$, the  formula $(\cal C \cal C^{\ell oc})$ is just 
  $$\sigma_d(X_0)=\Theta_d(X_0)=\Lambda_d^{\ell oc}(X_0) .$$ 
 
 In the case    $X$ is a complex analytic subset of $\CC^n$, $\cal G=\widetilde G(d/2,n)$ (the $d/2$-dimensional complex vector subspaces of $\CC^n$) and $G=U_n(\CC)$, since by definition the number of points in the fibre of a projection of the germ $X_0$ onto a generic $d/2$-dimensional complex vector subspace of $\CC^n$ is the local multiplicity $e(X,0)$ of $X_0$, formula $(\cal C \cal C^{\ell oc})$  gives
$$e(X,0)=\sigma^{\tilde G(d/2,n)}_d(X_0)=\Theta_{d}(X_0).$$ 
Obtaining the equality $e(X,0)=\Theta_d(X_0)$
 as a by-product of the formula  $(\cal C \cal C^{\ell oc})$ provides a new proof of  Draper's result (see \cite{Dra}).

\begin{remark}
When $(X^j)_{j\in \{0, \cdots, k\}}$ is a Whitney stratification of the closed set $X$ (see for instance \cite{Tro1} and \cite{Tro2} 
for a survey on regularity conditions for stratifications) and  $0\in X^0 $, $\sigma_i(X_0)=1$, for $i\le \dim( X^0 )$
(see \cite{CoMe}, Remark 2.9).
Therefore, to sum up, when $(X^j)_{j\in \{0, \cdots, k\}}$ is a Whitney stratification of the closed definable set $X$ and  $d_0$ is the dimension of the stratum containing $0$, one has $$\sigma_*(X_0)=(1, \cdots, 1, \sigma_{d_0+1}(X_0),
\cdots, \sigma_{d-1}(X_0), \Lambda_d^{\ell oc}(X_0)(X_0)=\Theta_d(X_0), 0, \cdots, 0).$$
\end{remark}

On figure 5 we represent the data taken into account in the computation of the invariant $ \sigma_i(X_0)$. 
Here, contrary to the computation of the $\Lambda_i^{\ell oc}(X_0)$ where all the domains $K^{P,\varrho}_\ell$ matter (see figure 4),
only the domains $K^{P,\varrho}_\ell$ having the origin in their adherence
 (these domains are coloured in red in figure 5) are considered, since only these domains appear as  $\chi(\pi_P^{-1}(y)\cap Z\cap 
\bar B_{(0,\varrho)})\cdot
\Theta_i((K_\ell^{P,\varrho})_0)$ in the computation of $\theta(\pi_{P_0*}
({\bf 1}_{X_0}))$. 

In particular the domains $K^{P,\varrho}_\ell$  
(in green on figure 5) defined by the critical values of the projection  $\pi_P$ restricted to the link $X\cap S_{(0,\varrho)}$
are not considered in the definition of $\theta(\pi_{P_0*}
({\bf 1}_{X_0}))$. One can indeed  prove (see \cite{CoMe}, Proposition 2.5) that for any generic projection   $\pi_P$ exists $r_P>0$, such that  for all
$\varrho$, $0<\varrho<r_P$, 
the discriminant of the restriction of $\pi_P$ to the link $X\cap S_{(0,\varrho)}$ is at  a positive distance from $0=\pi_P(0)$.

\eject
$$\hskip-9,3cm  X\cap \bar B^n_{(0,\varrho)}$$
\vskip8.8cm
\vglue0cm
\hskip0cm\includegraphics{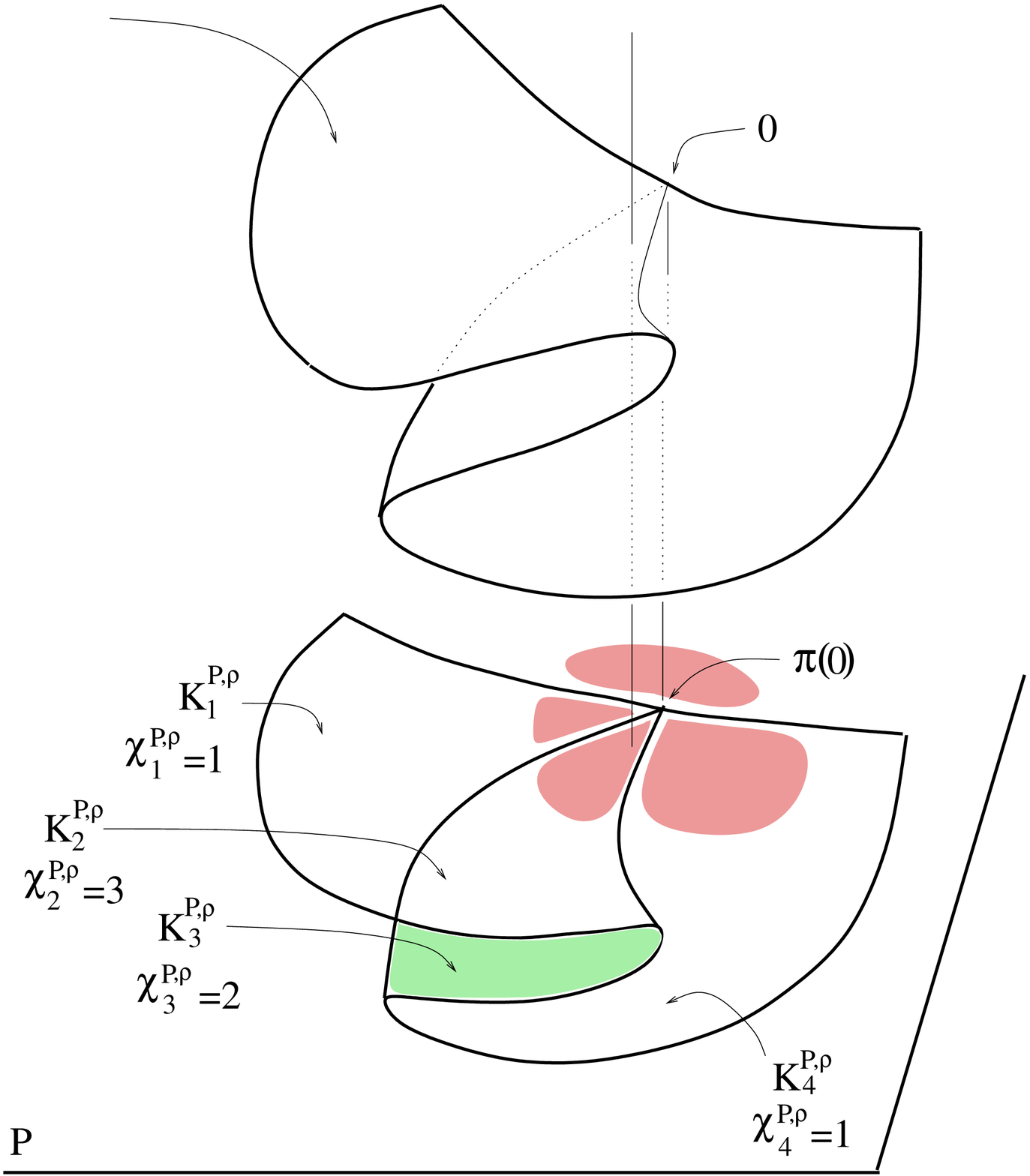}
\vglue0cm
\vskip0,5cm
\centerline{\bf fig. 5}
 \vskip0,8cm

Let us now deal with the question of what kind of invariants
of complex singularities the sequence $\sigma_*$ of invariants of real singularities  generalizes. 
For this goal, we consider that $X$ is a complex analytic 
subset of $\CC^n$ of complex dimension  $d$. 
One may define, like in the real case, the polar invariants of the germ $X_0$, denoted   
 $\tilde \sigma_i$, $i=0, \cdots, n$. These invariants are defined by generic projections 
on  $i$-dimensional complex vector subspaces of $\CC^n$. 
In the complex case, assuming  $0\in X$, there exists $r>0$, 
such that for $y$ generic in a generic $i$-dimensional vector 
space $P$ of $\CC^n$ and $y$ sufficiently closed to $0$
$$ \tilde \sigma_i(X_0)=\chi(\pi_P^{-1}(y)\cap X\cap \bar B_{(0,r)}).$$
In particular, as already observed,   $\tilde \sigma_d(X_0)$
is $e(X,0)$, the local multiplicity of $X$ at $0$.

In the case where $X$  is a complex hypersurface $f^{-1}(0)$, given by an analytic function $f:(\CC^n,0) \to (\CC,0)$ having at $0$
an isolated singularity, one has for a generic $y$ in the germ at $0$ of a generic $i$-dimensional vector space $P_0$ of $\CC^n$ 
$$ \chi (\pi_P^{-1}(y)\cap X\cap \bar B_{(0,\eta)})= \chi
(\pi_P^{-1}(0)\cap f^{-1}(\varepsilon)\cap \bar B_{(0,\eta)}),   $$
where $\varepsilon$ is generic in $\CC$, sufficiently close to
$0$ and  $0<\vert \varepsilon\vert \ll   \eta  \ll 1$.
Therefore, for   $0<\vert \varepsilon\vert \ll\eta\ll 1$, the integer $ \chi(\pi_P^{-1}(0)\cap f^{-1}(\varepsilon)\cap
\bar B_{(0,\eta)})   $
is the Euler-Poincar\'e characteristic of the Milnor fibre of $f$ restricted to  $P^\perp$, that is to say $1+(-1)^{n-i-1}\mu^{(n-i)} $.  In the case where $X$ is a complex analytic hypersurface
of $\CC^n$ with an isolated singularity at  $0$, we thus have
$$\widetilde{\sigma_i}(X_0)=1+(-1)^{n-i-1}\mu^{(n-i)}.$$

For $X$ a complex analytic subset of $\CC^n$ of dimension $d$, 
$d$ being  not necessarily $n-1$, the complex invariants
$\widetilde{\sigma}_i(X_0)$
have been first considered by  Kashiwara in  \cite{Ka}
(where the balls are open and not closed as it is the case here). 
 An invariant $E_{X_0}^0$ is then defined in \cite{Ka} by induction on the dimension of  $X_0$ using $\tilde \sigma_i$. This invariant is studied in  \cite{Dub1}, \cite{Dub2} and in \cite{BrDuKa} where a multidimensional version  $E^k_{X_0}$ of  $E_{X_0}^0$ is given (see also \cite{Me}). The definition is the following

$$ E_{X_0}^k = \sum_{X^{j_0}\subset \bar X^j\setminus X^j, 
\ \dim(X^j)<\dim(X_0)} E^k_{\bar X^j}\cdot
\tilde \sigma_{k+\dim(X^j)+1}(X_0),$$
where $(X^j)$ is a Whitney stratification of $X_0$, $X^{j_0}$ the  stratum containing $0$ and $E^k_{\{0\}}=1$.
The authors then remark that   (see also \cite{Dub1}, \cite{Dub2})
$$ E_{X_0}^k = Eu^k_{X_0},$$
where 
 $Eu^k_{X_0}=Eu_{(X_0\cap H)}$,  $H$ is a general vector subspace of dimension  $k$ of $\CC^n$ and where $Eu$ is the 
 local Euler obstruction of $X$ at $0$,
introduced by  MacPherson in  \cite{MacPh}. In particular, 
 $$E^0_{X_0}=Eu_{X_0}.$$
 
Let us now denote $\cal P^i(X_0)$, $i=0, \cdots, d$, the codimension $i$ polar variety of $X_0$, that is to say the closure
 of the critical locus of the projection of the regular part of
$X_0$ to a generic vector space of $\CC^n$ of dimension $d-i+1$.  
 The following relation between invariants $E^k_{X_0}$ and the local multiplicity of the polar varieties $\cal P^i(X_0)$ 
 is obtained in \cite{BrDuKa}
$$ (-1)^i(E_{X_0}^{\dim(X_0)-i-1}-E_{X_0}^{\dim(X_0)-i})= 
e(\cal P^i(X_0),0), $$
which in turn gives (see also \cite{Me}, \cite{LeTei1}, \cite{LeTei2}, \cite{LeTei3}, \cite{Dub2}) 
$$Eu_{X_0}=\sum_{i=0}^{d-1}(-1)^i e(\cal P^i(X_0),0), $$  
where $e(\cal P^i(X_0),0)$ is, as before, the local  multiplicity at  $0$ of the codimension $i$ polar variety  
$\cal P^i(X_0)$ of  $X_0$ at $0$.

All the invariants $\tilde \sigma_i(X_y)$, $E_{X_y}^i $, $e(\cal P^i(X_y),y)$, viewed as functions of the base-point $y$, enjoy the same remarkable property: they 
can  detect subtle variations of the geometry of an analytic family $(X_y)$, in the sense that the family  $(X_y)$ may be  Whitney stratified with $y$ staying in the same stratum if and only if these invariants are constant with respect to the parameter $y$. 
Without proof, it is actually stated in 
\cite{Dub1} Proposition 1, \cite{Dub2} Theorem II.2.7 page 30,
 and  \cite{BrDuKa}, that the invariants $\tilde \sigma_i(X_y)$ are constant as
$y$ varies in a stratum of a Whitney stratification of $X_0$
(see also \cite{CoMe} Corollary 4.5). And in \cite{HeMe1}, \cite{Na2}, \cite{Tei2}
it is proved that the constancy of the multiplicities
$e(\cal P^i(X_y),y)$ as $y$ varies in a stratum of a stratification of $X_0$, is equivalent to 
the Whitney regularity of this stratification, giving also a proof,
considering the relations  between $e(\cal P^i(X_y),y)$ and $\tilde \sigma_i(X_y)$ stated above,
of the constancy of  $y\ni\to \tilde \sigma_i(X_y)$ along Whitney strata.

We sum-up these results in the following theorem, where
$e( \Delta^i(X_y),y)$ is the local multiplicity at $y$ of the discriminant
$\Delta^i(X_y)$ associated to $\cal P^i(X_y)$, that is the image of $\cal P^i(X_y)$ under the generic projection that gives rise to $\cal P^i(X_y)$.
 
\begin{thm}[\cite{HeMe1}, \cite{Na2}, \cite{LeTei3}, \cite{Tei2} ]\label{Equising}
Let  $X_0$ be a complex analytic germ at $0$ of $\CC^n$
endowed with a stratification $(X^j)$.
The following statements are equivalent  
 
\begin{enumerate}
\item The stratification  $(X^j)$ is a  Whitney stratification.
\item The functions $X^j\ni y\mapsto e(\cal P^i(X^\ell_y),y) $, for $i=0, \cdots, d-1$ and any pairs $(X^j,X^\ell)$ such that
$X^j\subset \overline{X^\ell}$, are constant.
\item The functions $X^j\ni y\mapsto e( \Delta^i(X^\ell_y),y) $,  for $i=0, \cdots, d-1$ and for any pairs $(X^j,X^\ell)$ such that
$X^j\subset \overline{X^\ell}$, are constant.
\item The functions $y\ni X^j \mapsto \tilde \sigma_i(X_y) $,  for $i=1, \cdots, d$  are constant.
\end{enumerate} 

\end{thm}

\begin{remark}
In the real case the functions $y  \mapsto  \sigma_i(X_y) $ are not $\ZZ$-valued functions as in the complex case, but   $\RR$-valued functions 
and 
in general one can not stratify a compact definable set in such a way that the restriction of these functions to the strata are constant. 
However, it is proved in \cite{CoMe} Theorems 4.9 and 4.10, that Verdier regularity for a stratification implies continuity of the restriction
of $y  \mapsto  \sigma_i(X_y) $ to the strata of this stratification. Since, in the complex setting, Verdier regularity is the same as Whitney regularity, 
this result is the real counterpart of Theorem \ref{Equising}. Note that in the real case one can not expect that the continuity or even the constancy 
of the functions $y  \mapsto  \sigma_i(X_y) $ in restriction to the strata of a given stratification  implies a convenient regularity condition for this stratification (see the introduction of \cite{CoMe}). 
\end{remark}

As a conclusion of this section, the complex version $\tilde \sigma_i$ of the real polar invariants 
$\sigma_i$ of definable singularities plays a central role in singularity theory since they let us compute classical invariants
of singularities and since their constancy, with respect to the parameter of an analytic family, as well as the constancy of other  classical invariants related to them, means that the family does not change its geometry. 
Since our polar invariants  $ \sigma_i$ appear now as the real counterpart of classical complex invariants, we'd like to understand in the sequel how they  are related to the local Lipschitz-Killing invariants $\Lambda_i^{\ell oc}(X_0)$ coming from the differential geometry of the deformation of the germ $X_0$ through its tubular neighbourhoods family. 

This is the goal of the next section.

\subsection{ Multidimensional local Cauchy-Crofton formula }\label{section 2.3}
The local Cauchy-Crofton formula $(\cal C \cal C^{\ell oc})$ given at  \ref{LCCF} already equals $\sigma_d$ and 
$\Lambda_d^{\ell oc}$ over definable germs. This relation suggests
a more general relation  between the  $\Lambda^{\ell oc}_i$'s and the  $\sigma_j$'s. We actually can prove that each invariant of one family is a linear combination 
of the invariants of the other family. The precise statement is given by the following formula $(\cal C\cal C^{\ell oc}_{mult})$.  
 
\begin{MLCCF}[\cite{CoMe} Theorem 3.1]  
There exist real numbers
$(m_i^j)_{1\le i,j\le n,i<j }$ such that, for any definable germ 
 $X_0$, one has 

$$
\begin{pmatrix}
\Lambda^{\ell oc}_1(X_0) \\ \vdots \\ \Lambda^{\ell oc}_n(X_0)
\end{pmatrix}
=
\begin{pmatrix} 1 & m_1^2  & \ldots & m_1^{n-1}  & m_1^n \\
                0 & 1      & \ldots & m_2^{n-1 } & m_2^n  \\
           \vdots &        &        &            & \vdots \\
                0 & 0      & \ldots &           0& 1 \\ 
\end{pmatrix}
                     \cdot
\begin{pmatrix} \sigma_1(X_0) \\ \vdots \\ \sigma_n(X_0) \\
\end{pmatrix} \eqno (\cal C\cal C^{\ell oc}_{mult})$$
\vskip2mm
\noindent 
These constant real numbers are given by
 $   \displaystyle  m_i^j=
  {\alpha_j  \over \alpha_{j-i}\cdot
\alpha_i} 
\begin{pmatrix} i\\ j\\
\end{pmatrix} -   {\alpha_{j-1}  
\over \alpha_{j-1-i}\cdot
\alpha_i}
\begin{pmatrix} i \\ j-1\\
\end{pmatrix},$
for $   i+1\le j \le n.$
\end{MLCCF}

\begin{remark}
Applied to a $d$-dimensional definable 
germ $X_0$, the last {\sl a priori} non-trivial equality provided by   formula 
$(\cal C\cal C^{\ell oc}_{m ult})$, involving the $d$-th line of the matrices, is 
$$\Theta_d(X_0)=\Lambda_d^{\ell oc}(X_0)= \sigma_d(X_0),$$
 which is the local 
Cauchy-Crofton formula $(\cal C\cal C^{\ell oc})$.  
The local Cauchy-Crofton formula $(\cal C\cal C^{\ell oc})$ expresses
the $d$-density of a $d$-dimensional germ as the mean value of the number of points 
in the intersection of this germ with a $(n-d)$-dimensional affine space of $\RR^n$.
This number of points may be viewed as the Euler-Poincar\'e characteristic of this intersection. 
Now, since for $d$-dimensional germs  the $d$-density is the last 
invariant of the sequence $\Lambda_*^{\ell oc}$  and since
formula $(\cal C\cal C^{\ell oc}_{m ult})$ expresses all invariants $\Lambda_i^{\ell oc}$  in terms of the mean values of 
the Euler-Poincar\'e characteristics of the multidimensional plane sections of our germ,   we see in
formula $(\cal C\cal C^{\ell oc}_{m ult})$ a multimensional version of the local Cauchy-Crofton formula.

\end{remark}

\subsection{Valuations theory and Hadwiger principle}\label{section 2.4}
In the previous section \ref{section 2.3}, with formula $(\cal C \cal C^{\ell oc}_{mult})$, we have answered 
the question: {\sl how are the local Lipschitz-Killing invariants and the  polar invariants  related?} 
In this section  we would like to risk some speculative and maybe prospective insights about 
formulas similar to  $(\cal C\cal C^{\ell oc}_{m ult})$, that is to say formulas 
 linearly relating two families of local invariants of singularities. We include in this scope 
 the complex formulas presented in Section \ref{section 2.1}. 
 For this goal we first recall some definitions and celebrated statements from convex geometry, since it appears 
 that from the theory of valuations on convex bodies one can draw precious lessons on the question: {\sl  why our additive invariants are linearly dependent?}

We have already observed that the additive functions $\Lambda_i$, $\Lambda_i^{\ell oc}$ and $\sigma_j$
are invariant under isometries of $\RR^n$. Therefore the first question we would like to address to convex geometry is the following: to what extend those invariants are models of additive and rigid motion invariants? 

 The systematic study of additive invariants (of compact convex sets of $\RR^n$) has been inaugurated by 
 Hadwiger and his school and motivated by Hilbert's third problem (solved by Dehn by introducing the so-called Dehn invariants) consisting in classifying scissors invariants of polytopes (see for instance \cite{Car} for a quick introduction to Hilbert's third problem). 
One of the most striking results
in this field is  Hadwiger's theorem that characterises the set of additive and 
rigid motion invariant functions (on the set of compact convex subsets of $\RR^n$) as
the vector space spanned by the 
 $\Lambda_i$'s. We give now the needful definitions to state Hadwiger's theorem and the still-open 
 question of its extension to the spherical case (for more details one can refer to \cite{McMu-Sch}, \cite{Sch3} or \cite{Sch4}).

We denote by  $\cal K^n$  (resp. $\cal K\SS^{n-1}$) the set of compact convex sets of   
$\RR^n$ (resp. of $\SS^{n-1}$, that is to say the intersection of the sphere $\SS^{n-1}$ and conic compact convex sets of   $\cal K^n$ with vertex the origin of $\RR^n$).
A function  $v: \cal K^n \to \RR$ (resp.
$v: \cal K\SS^{n-1}\to \RR $) is called a valuation (resp. a spherical valuation) when $v(\emptyset)=0$ and
for any  $K, L \in \cal K^n$
(resp. $K, L \in \cal K\SS^{n-1}$) such that  $K\cup L\in \cal K^n$
(resp. $K\cup L \in \cal K\SS^{n-1}$), one has the additivity property
$$ v(K \cup L)=v(K) + v(L) - v(K\cap L).  $$
One says that a valuation  $v$ on $\cal K^n$ (resp. $\cal K\SS^{n-1}$) is continuous
when it is continuous with respect to the Hausdorff metric on
$\cal K^n$ (resp. on $\cal K\SS^{n-1}$).
A valuation $v$ on $\cal K^n$ (resp. $\cal K\SS^{n-1}$) is called simple
when the restriction of 
 $v$ to convex sets with empty interior is zero.
Let  $G$ be a subgroup of the orthogonal group $O_n(\RR)$.
 A valuation   $v$ on $\cal K^n$ (resp. $\cal K\SS^{n-1}$) is  
{\sl $G$-invariante} when it is invariant under the action of translations of $\RR^n$ and the action of $G$
on  $\cal K^n$
(resp. the action of  $G$ on  $\cal K\SS^{n-1}$).

 The Hadwiger theorem emphasizes the central role played in convex geometry by the 
 Lipschitz-Killing invariants as additive rigid motion invariants. 
 
\begin{thm}[\cite{Had}, \cite{Kla}]\label{Had}
 A basis of the vector space of $SO_n(\RR)$-invariant and continuous valuations on  $\cal K^n$ is
$(\Lambda_0=1,\Lambda_1,\cdots, \Lambda_n=Vol_n)$.

Equivalently (by an easy induction argument) a basis of the vector space of continuous and $SO_n(\RR)$-invariant simple 
valuations on  $\cal K^n$ is $\Lambda_n=Vol_n$. 
\end{thm}

This statement forces a family of $n+2$ additive, continuous and $SO_n(\RR)$-invariant functions
on euclidean convex bodies to be linearly dependent.
The second formulation of Hadwiger's theorem, concerning the space of simple valuations,  
enables to address the question of such a rigid structure of the space of valuations in the
setting of spherical convex geometry. 
This question of a spherical version of Hadwiger's result has been address by  Gruber and Schneider. 

\begin{question}[\cite{Gr-Sc} Problem 74, \cite{McMu-Sch} Problem 14.3]\label{Hadquestion}
 Is a simple, continuous and $O_n(\RR)$-invariant valuation on $\cal K\SS^{n-1}$ 
 a multiple of the $(n-1)$-volume on $S^{n-1}$?
\end{question}

\begin{remark}
In the case $n\le 3$ a positive answer to this question given in \cite{McMu-Sch}, Theorem 14.4, and in the 
easy case where the simple valuation has constant sign one also has a positive answer given in \cite{Sch1} Theorem 6.2, 
and \cite{Sch2}. Note that in this last case the continuity is not required and that the valuation is {\sl a priori} 
defined only on convex spherical polytopes.  
\end{remark}

This difficult and still unsolved problem
naturally appears as soon as one consider the localizations $\Lambda^{\ell oc}_i$ of $\Lambda_i$ and their 
relation with other classical additive invariants such as the $\sigma_j$'s. Indeed, the question of why and how such invariants are related falls within the framework of Question \ref{Hadquestion}. Let's clarify this principled position.

The invariants $(\Lambda^{\ell oc}_i)_{i\in \{0, \cdots, n\}}$ define
spherical $O_n(\RR)$-invariant and continuous valuations $(\widehat\Lambda_i)_{i\in \{0, \cdots, n\}}$ 
on the convex sets of $\SS^{n-1}$ by the formula
$$ \widehat\Lambda_i(K):=\Lambda^{\ell oc}_i(\widehat K_0)={1\over 
\alpha_i}\Lambda_i(\widehat K\cap \bar B_{(0,1)}),\eqno(\widehat 8)$$
where $K$ is a convex set of $\SS^{n-1}$, that is to say the trace in 
  $\SS^{n-1}$ of the cone $\widehat K=\RR_+\cdot K$ with vertex the origin of $\RR^n$.
  Another possible finite sequence of continuous and $O_n(\RR)$-invariant spherical valuations 
on  convex polytopes of  $\SS^{n-1}$ is
$$\Xi_i(P):=\sum_{F\in {\cal F}_i(P)} 
Vol_i(F)\cdot\gamma(\widehat F,\widehat P)=Vol_i(S^i(0,1))
\sum_{F\in {\cal F}_i(P)} 
\Theta_i(\widehat F_0)\cdot\gamma(\widehat F,\widehat P)
, \eqno(\widehat 9)$$
where $P\subset S^{n-1}$ is a spherical  polytope, that is to say that 
$\RR_+\cdot P=\widehat P$ is the intersection of a finite number of closed  half vector spaces 
of $\RR^n$, $\cal F_i(P)$ the set of  all $i$-dimensional faces of  $P$ (the
$(i+1)$-dimensional faces of $\widehat P$) and
$\gamma(\widehat F,\widehat P)$ the external angle of $\widehat P$
along $\widehat F$. The valuations $\Xi_i$ are the natural spherical substitutes of the euclidean Lipschitz-Killing 
curvatures $\Lambda_i$ according to formula $(2)$.

Finally the polar invariants $(\sigma_i)_{i\in \{0, \cdots, n\}}$
also define continuous and  $O_n(\RR)$-inva\-riants  spherical valuations  
$(\widehat{\sigma_i})_{i\in \{0, \cdots, n\}}$ on the convex sets of 
$\SS^{n-1}$, according to a formula of the same type that formula 
$(\widehat8)$ 
$$ \widehat{\sigma_i}(K):=\sigma_i(\widehat K_0). \eqno(\widehat{10}) $$
These three families of continuous and $O_n(\RR)$-invariant spherical  valuations  
 $$(\widehat\Lambda_i)_{i\in \{0, \cdots, n\}},
({\Xi_i})_{i\in \{0, \cdots, n\}},
(\widehat{\sigma_i})_{i\in \{0, \cdots, n\}}$$ being linearly independent families
in the space of spherical  valuations, a  positive answer to Question \ref{Hadquestion} 
would have for direct consequence  that each element of one family is a linear 
combination of elements of any of the other two families. 

Therefore, in restriction to polyhedral cones, each element of the family   
$(\Lambda_i^{\ell oc})_{i\in \{0, \cdots, n\}}$ could be expressed as a linear combination  
(with universal coefficients)
of elements of the family $(\sigma_i)_{i\in \{0, \cdots, n\}}$
and conversely. Despite the absence of any positive answer to Question \ref{Hadquestion} for $n>3$, this linear dependence is proved in  \cite{CoMe} (Theorem A4 and A5) over the set of convex polytopes (and 
in \cite{CoMe}, section 3.1, even over the set of definable cones). It is actually shown that 
each invariant  $\Lambda^{\ell oc}_i$ and each invariant 
$\sigma_j$ may be expressed as a linear combination (with universal coefficients) of elements of the family $({\Xi_i})_{i\in \{0, \cdots, n\}}$. The coefficients involved in such linear combinations may be explicitly computed by considering the case of polytopes.

 In conclusion, an anticipating positive answer to this question would imply the existence of a finite  number 
of independent models for additive, continuous  and
$O_n(\RR)$-invariant functions on convex and conic germs. Consequently, solely following this principle and
restricted at least to convex cones, our local invariants $\sigma_j$ and $\Lambda^{\ell oc}_i$ would be automatically linearly dependent. In the next 
step, in order to extend a relation formula involving some valuations from the set of  finite union of convex conic polytopes to the set of general 
conic definable set one has to prove some general statement according to which the normal cycle of a definable conic set may be approximate by the 
normal cycles of a family of  finite union of convex conic polytopes. We do not want to go more into technical details and even define the notion of
 normal cycle introduced by  Fu; we just point out that this issue has been tackled recently in \cite{FuSco}. 
 Finally to extend a relation formula from the set of conic definable set to the set of all definable germs one just has to use the local conic 
 structure of definable germs and the deformation on the tangent cone (see \cite{CoMe}, \cite{Du12}).

To finish to shed light on convex geometry as an area from which some strong relations between singularity invariants may be understood, let us remark that the following generalization of Hadwiger's theorem 
\ref{Had} has been obtained by  Alesker.

%

\begin{thm}
 Let  $G$ be a compact subgroup of $O_n(\RR)$. 
\begin{enumerate}

\item The vector space 
$Val^G(\cal K^n)$ of continuous, translation and $G$-invariant valuations on $\cal K^n$ has 
finite dimension if and only if  
$G$ acts transitively on $\SS^{n-1}$ (see \cite{Ale1} Theorem 8.1,
\cite{Ale5} Proposition 2.6).

\item One can endowed the vector space  $Val^G(\cal K^n)$ with a product 
(see \cite{Ale4}, \cite{Ale-Fu}) providing a graded algebra structure
 (the graduation coming from the homogeneity degree of the valuations) and 
$$ \begin{matrix} \RR[x]/(x^{n+1}) & \to & Val^{O_n(\RR)}(\cal K^n)
=Val^{SO_n(\RR)}(\cal K^n)\\
                    \hfill x & \mapsto & \Lambda_1 \hfill\\
                    \end{matrix}  $$
is an   isomorphism of graded algebras
(see \cite{Ale4}, Theorem 2.6).

\end{enumerate}

\end{thm}
%

\section{Generating additive invariants via generating functions}\label{Section Mot}

The possibility of generating invariants from a deformation of a singular set into a family of approximating and less complicated sets is perfectly illustrated by the work developed by  Denef and  Loeser consisting in 
stating that some generating series attached to a singular  germ are rational. Such generating series have their coefficients in some convenient  ring reflecting the special properties of the invariants to highlight, such as additivity, multiplicativity, analytic invariance, the relation with some specific group action and so forth, and on the other hand each of these coefficients is attached to a single element of the deformation family. It follows that such a  generating series
captures the geometrical aspect that one aims to focus on through  the deformation family  as well as its rationality   
indicates that asymptotically this geometry  specializes on the geometry of the special fibre approximated by the deformation family. Indeed, being rational strongly expresses that a series is encoded by a finite amount of data concentrated in its higher  coefficients.   
To be more explicit we now roughly describe how   Denef and  Loeser define the notion of motivic Milnor fibre (for far more complete and precise introductions to motivic integration which is the central tool of the theory, and to motivic invariants in general, the reader may refer to \cite{Bli}, \cite{Clu-Loe2008}, \cite{Clu-Loe2011}, \cite{DenBour}, \cite{DenLoe1999}, \cite{DenLoe2001a}, \cite{DenLoe2001b}, \cite{DenLoe2002}, \cite{GorYaf}, \cite{GSLuMe}, \cite{Hal1}, \cite{Hal2}, \cite{Loe2009}, \cite{Loo}, \cite{Veys}). 
\subsection{The complex case}
A possible starting point of the theory of motivic invariants may be attributed to the works of Igusa (see \cite{Igu1}, \cite{Igu2}, \cite{Igu3})
on zeta functions introduced by Weil (see \cite{Wei}). In the works of Igusa the rationality of some generating Poincar\'e series is proved. These series, called zeta functions, have for coefficients the number of points in $\cal O/{\textswab M}^m$ of $f= 0 \mod {\textswab M}^m$, for $f$ a $n$-ary polynomial with coefficients in the valuation ring $\cal O$ of some discrete valuation field of characteristic zero, with maximal ideal ${\textswab M}$ and finite residue field of cardinal $q$. This result amounts to prove the rationality (as a function of $q^{-s}$, $s\in \CC$, $\Re e(s)>0$) of an integral of type
$$ \int_{\cal O^m} \vert f\vert^s \vert ds \vert\eqno{(11)} $$
(the   Igusa local zeta function) which is achieve using a convenient resolution of singularities of $f$, as described in the introduction (see also \cite{Den84} for comparable statements on Serre's series and a strategy based on Macintyre's proof of elimination of quantifiers as an alternative to Hironaka's resolution of singularities). 
In the opposite direction, the rational function, expressed in terms of the data of a resolution of $f$, associated to a polynomial germ $f:(\CC,0)\to (\CC,0)$ by the expression provided by the computation of Igusa's integral in the discrete valuation field case let 
Denef and Loeser define intrinsic invariants attached to the complex germ $(f^{-1}(0),0)$, called
topological zeta functions (see \cite{DenLoe1992}). Another key milestone in the systematic use of discrete valuation fields (here with finite residue field and more explicitely in the $p$-adic context) have been reached in Batyrev's paper \cite{Bat}, where it is shown that two birationally equivalent Calabi-Yau manifolds over $\CC$ have the same Betti numbers. Indeed, by Weil's conjectures, these Betti numbers are obtained from the rational expression of the local zeta functions having for coefficients the number of points in the reductions
modulo $p^m$ of the manifolds into consideration (viewed as defined over $\QQ_p$ when they are defined over $\QQ\subset \CC$) and on the other hand, these local zeta functions may be computed by Igusa's integrals over these manifolds. Being birationally equivalent, these manifolds provide the same integrals. 
Kontsevich,  in his seminal talk \cite{Kon},   extended this method (consisting in shifting a complex geometric problem 
in a discrete valuation field setting) in the equicharacteristic setting by developing an integration theory  
in particular over $\CC[[t]]$. Note that the theory may be developed in great generality and not only in equicharacteristic zero (see \cite{Loo}, \cite{Clu-Loe2008}).  
The idea of Kontsevich was to define an integration theory over arc spaces, say $\CC[[t]]$, by considering a measure with values in the Grothendieck ring $K_0(\Var_\CC)$  of algebraic varieties over $\CC$ (localized
by the multiplicative set generated by the class $\L$ of $\AA^1$ in $K_0(\Var_\CC)$). The main tool in this context 
being a change of variables formula that allows computation of integrals through morphisms, and in particular through a morphism  given by a resolution of singularities.  
Formally the ring  $K_0(\Var_\CC)$ is the
free abelian group generated by isomorphism classes $[X]$ of varieties $X$ over $\CC$, with the relations
$$ [X\setminus Y]=[X]-[Y],$$
for $Y$ closed in $X$, the product of ring being given by the product of varieties (see for instance \cite{NicSeb}). Denoting $\L$ the class of $\AA^1$ in $K_0(\Var_\CC)$,
we then denote $\cal M_\CC$ the localization $K_0(\Var_\CC)[\L^{-1}]$. 
Any additive and multiplicative invariant on $\Var_\CC$ with non zero value at $\AA^1$, such as the Euler-Poincar\'e characteristic or the Hodge 
characteristic (both with compact support), factorizes through the universal additive and multiplicative class 
$\Var_\CC\ni X\mapsto [X]\in \cal M_\CC$.  
Now we equip the space $\cal L(\CC^n,0)$ of formal arcs of $\CC^n$ passing through $0$ at $0$  with the 
above-mentioned measure that provides a $\sigma$-additive measure, with values in a completion $\hat{\cal M}_\CC$ of 
$\cal M_\CC,$ for sets of the boolean algebra of the so-called constructible sets of $\cal L(\CC^n,0)$. 
Finally denoting $\cal L_m(\CC^n,0)$, $m\ge 0$, the set of polynomial arcs of $\CC^n$ of degree $\le m$, passing 
through $0$ at $0$, and for $f:(\AA^n,0)\to (\AA^1,0)$ a morphism having a (isolated) singularity at the origin, 
inducing the morphism $f_m:\cal L_m(\CC^n,0)\to \cal L(\CC,0)$, we denote 
$$\cal X_{m,0,1}:=\{ \varphi\in \cal L_m(\CC^n,0); f_m\circ \varphi = t^m+\hbox{high order terms} \},$$
and we define (see for instance \cite{DenLoe2002}) the  motivic zeta function of $f$ by
$$Z_f(T):=\sum_{m\ge 1} [\cal X_{m,0,1}]\ \L^{-mn}T^m.$$
 This generating series appears as the $\CC[[t]]$-substitute of 
$p$-adic zeta functions introduced by Weil, and whose rationality, following Igusa, 
amounts to compute an integral of type $(11)$.   Inspired by the analogy of $Z_f(T)$ with Igusa integrals, 
and also using a resolution of singularities of $f$ as presented in 
the introduction and using the Kontsevich change of variables formula for this resolution applied to the 
coefficients of $Z_f$ viewed as measures of constructible sets, Denef and Loeser proved (see \cite{DenLoe1998}, 
\cite{DenLoe1999}, \cite{DenLoe2002}) the rationality of $Z_f(T)$. 
With the notation given in the introduction, we then have 
$$Z_f(T)= \sum_{I\cap \cal K \not= \emptyset}(\L-1)^{\vert I\vert -1 }
[\widetilde E^0_I ]\prod_{i\in I}
\frac{\L^{-\nu_i}T^{N_i}}{1-\L^{-\nu_i}T^{N_i}} \eqno{(12)}$$
 where $\widetilde E^0_I$ is a covering of $E_I^0$ defined in the following way. Let $U$ be some affine open subset of $M$ such that on $U$, $f\circ \sigma(x)=u(x)\prod_{i\in I}x_i^{N_i}$, with $u$ a unit. Then 
 $\widetilde E^0_I$ is obtained by gluing along $E_I^0\cap U$ the sets 
 $$\{ (x,z)\in (E_I^0\cap U)\times \AA^1; z^{m_I}\cdot u(x)=1\}, $$ 
 where $m_I=\hbox{gcd}(N_i)_{i\in I}$.

 \begin{question}[Monodromy conjecture]
 We do not know how the poles (in some sense) of the rational expression of $Z_f$ relates on the eigenvalues of the 
 monodromy function $M$ associated to the singular germ $f:(\CC^n,0)\to (\CC,0)$. The Monodromy Conjecture
 of Igusa, that has been stated in many different forms after Igusa, asserts that when $\L^\nu-T^N$ 
 indeed appears as denominator of the rational expression of $Z_f$ viewed as an element of the ring 
 generated by $\hat{\cal M}_\CC$ and $T^N/(\L^\nu-T^N)$, $\nu,N>0$, 
then  
 $e^{2i\pi \nu/N}$ is an eigenvalue of $M$  (see for instance \cite{DenBour} and the references given in this
 article for additional classical references).
 \end{question}

 \begin{remark}
 The rationality of $Z_f$ illustrates again a deformation principle; the family $(\cal X_{m,0,1})_{m\ge 1}$ may be considered as a family of tubular neighbourhoods in $\cal L (\CC^n,0)$ around the singular fibre $\cal X_0=\{f=0\}$ and in a neighbourhood of the origin, 
 with respect to the ultrametric distance given by the order of arcs. Now the rationality 
of $Z_f$ expresses the regularity of the degeneracy of the geometry of $\cal X_{m,0,1}$ onto the geometry of 
$\cal X_0$. Following this principle, the rational expression $(12)$ of
$Z_f$ is supposed to concentrate the part of the geometrical information encoded in $(\cal X_{m,0,1})_{m\ge 1}$ that accumulates at infinity in the series $Z_f$.
\end{remark}

 This is achieved in particular 
by the following observation (see \cite{DenLoe1998}, \cite{DenLoe2001b}, \cite{DenLoe2002}): the negative 
of the   
constant term of the formal expansion as a power series in $1/T$
of the rational 
expression of $Z_f$ given by formula $(12)$
defines the following element in $\hat{\cal M}_\CC$
$$ S_f:=\sum_{I\cap \cal K \not= \emptyset}(\L-1)^{\vert I\vert -1 } \
[\widetilde E^0_I ],$$ 
called the motivic Milnor fibre of $f=0$ at the singular point $0$ of $f$. Taking the realization of 
$S_f$ under the morphism $\chi:  \hat{\cal M}_\CC \to \ZZ$
(note that $\chi(\L)=1$) gives, in particular, by the A'Campo formula recalled in the introduction,
$$ \chi(S_f)=
\sum_{ i\in \cal K} N_i \cdot \chi(E_{\{i\}}^0)=\chi(X_0)=1+(-1)^{n-1}\mu.\eqno(13)$$

\begin{remark}  Generally speaking, taking the constant term in the expansion of a rational function $Z$ as a power series in $1/T$
, amounts to consider $\lim_{T\to \infty }Z(T)$ (in a setting where this makes sense). 
A process that gives an increasing importance to the $m$-th coefficient of $Z$ as $m$ itself increases.
On the other hand, the coefficients of $Z_f$ may be directly interpreted
at the level of the Euler-Poincar\'e characteristic, which could be seen as the first topological degree
of realization of $K_0(\Var_\CC)$, and it turns out that the sequence $(\chi(\cal X_{m,0,1}))_{m\ge 1}$ has a strong regularity
since it is in fact 
periodic. Indeed, one has by \cite{DenLoe2002} Theorem 1.1
$$\chi(\cal X_{m,0,1})=\Lambda(M^m), \ \forall m\ge 1 \eqno(14)$$  
and by quasi-unipotence of $M$ 
(see \cite{SGA7} I.1.2) there exists $N>1$ such that
the order of the eigenvalues of $M$ divides $N$. It follows from $(14)$ that 
$\chi(\cal X_{m+N,0,1})=\chi(\cal X_{m,0,1}), m\ge 1$. 
\end{remark}    
   
Now formula $(13)$, showing that the motivic Milnor fibre has 
a realization, via the Euler-Poincar\'e characteristic, on the Euler-Poincar\'e characteristic of the set-theoretic Milnor fibre, is a direct 
consequence of formulas $(14)$ (note in fact that the proofs of $(13)$ and $(14)$, using  A'Campo's formulas and a resolution of singularities, are 
essentially the same and thus gives comparable statements). Indeed, as noticed by Loeser (personal communication), working with $\chi$ 
instead of formal classes
of $K_0(\Var_\CC)$, on gets by definition of $Z_f$  
$$\chi(Z_f)=\sum_{m\ge 1}\chi(\cal X_{m,0,1})T^m$$ that gives in turn, 
 by formula $(14)$, $$\chi(Z_f)=\sum_{m\ge 1}\Lambda(M^m)T^m=
\sum_{m=1}^N\Lambda(M^m)\sum_{k\ge 0}T^{m+kN}=\sum_{m=1}^N\Lambda(M^m)\frac{T^m}{1-T^N}.$$
Since $\chi(S_f)=-\lim_{T\to \infty}\chi(Z_f)$, on finally find again that  
$$\chi(S_f)=\Lambda(M^N)=\Lambda(Id)=\chi(X_0)=\chi(\bar X_0)=1+(-1)^{n-1}\mu.$$
   
 \begin{remark}
 One may consider a more specific Grothendieck ring, that is to say a ring with more relations, in order to take into account the monodromy action on the Milnor fibre. In this equivariant and more pertinent ring   equalities 
 $(12)$ and $(14)$ are still true (see \cite{DenLoe2001b}, \cite{DenLoe2002} Section 2.9)
 \end{remark}   
   
\begin{remark}  In \cite{HruLoe}, Hrushovski and Loeser gave a proof of equality $(14)$ without using a resolution of singularity, and therefore without using A'Campo's formulas. Since a computation of $Z_f$ in terms of the data associated to a particular resolution of the singularities of $f$ leads to the simple observation that one computes in this way an expression already provided by A'Campo's formulas, the original proof of $(14)$ may, in some sense, appear as a not direct proof.  The proof proposed in \cite{HruLoe} uses \'etale cohomology of non-archimedean spaces and motivic integration in the model theoretic version of \cite{HruKaz1} and \cite{HruKaz2}. 
\end{remark}

\begin{remark} To finish with the complex case, let us note that in \cite{Sai} and \cite{SteZuc} a mixed Hodge structure on 
the Milnor fibre $f^{-1}(t)$ at infinity ($\vert t \vert \gg 1$) has been defined by a deformation process, letting $t$ goes 
to infinity (see also \cite{Sab}). In \cite{RaiCRAS} and \cite{RaiSMF} a corresponding motivic Milnor fibre $S_f^\infty$ has then 
be defined.
\end{remark}

\subsection{The real case}
A real version of $(12)$, giving rise to a real version of $(13)$ has been obtained
in \cite{CoFi} (see also \cite{Fi}). In the real case, a singular germ $$f:(\RR^n,0)\to (\RR,0) $$ defines  two smooth bundles
$(f^{-1}(\varepsilon)\cap B_{(0,\eta)})_{ 0<-\varepsilon\ll \eta\ll 1} $ and $(f^{-1}(\varepsilon)\cap B_{(0,\eta)})_{ 0<\varepsilon\ll \eta\ll 1}$
and as well as $\cal X_{m,0,-1}$ and  $\cal X_{m,0,1}$, it is natural to consider the two sets 
$$\cal X_{m,0,>}:=\{ \varphi\in \cal L_m(\RR^n,0); f_m\circ \varphi = at^m+\hbox{high order terms}, a>0 \},$$
 and 
$$\cal X_{m,0,<}:=\{ \varphi\in \cal L_m(\RR^n,0); f_m\circ \varphi = at^m+\hbox{high order terms}, a<0 \}.$$
Let us denote $X_0^{-1}$ and $X_0^{+1}$ the fibre $f^{-1}(\varepsilon)\cap B_{(0,\eta)}$ for respectively $\varepsilon<0$ and $\varepsilon >0$. 
 While $\cal X_{m,0,a}$, for $a\in \CC^\times$, is a constructible set having 
a class in $K_0(\Var_\CC)$, the sets $\cal X_{m,0,<}$ ans $\cal X_{m,0,>}$ are real semialgebraic sets and unfortunately, the Grothendieck ring of real semialgebraic sets is the trivial ring $\ZZ$, since semialgebraic sets admit semialgebraic cells decomposition.  

Therefore, in the real case, since we have to deal with two signed Milnor fibres, we cannot mimic the construction of $K_0(\Var_\CC)$.  To overcome this issue, in \cite{Fi} we proposed to work in the Grothendieck ring of real (basic) semialgebraic formulas, $K_0(\hbox{BS}_\RR)$. In this ring no semialgebraic isomorphism relations between semialgebraic sets, but algebraic isomorphim relations between sets given by algebraic formulas,  are 
imposed and distinct real semialgebraic formulas 
having the same set of real points in $\RR^n$ may have different classes. In particular, a first order basic 
formula in the language of ordered rings with parameters from $\RR$ may have a nonzero class in $K_0(\hbox{BS}_\RR)$ whereas 
no real point satisfies it. The ring $K_0(\hbox{BS}_\RR)$ may be sent to the more convenient ring 
$K_0(\Var_\RR)\otimes \Z[\frac{1}{2}]$, where explicit computations of classes of basic semialgebraic formula are possible as long as 
 computations of classes of real algebraic formulas in the classical Grothendieck ring of real algebraic varieties $(K_0(\Var_\RR)$ are possible.

In this setting, since $\cal X_{m,0,>}$ and 
$\cal X_{m,0,<}$ are given by explicit basic semialgebraic formulas, they do have natural classes in $K_0(\hbox{BS}_\RR)$ and this
 allows the consideration of the associated zeta series 
$$Z_f^{?}=\sum_{m\ge 1} [\cal X_{m,0,?}]\ \L^{-mn}T^m \in (K_0(\Var_\RR)\otimes \Z[\frac{1}{2}])[\L^{-1}][[T]],   \ ?\in \{-1,+1, <,>\}.$$
It is then proved, with the same strategy as in the complex case (using a resolution of singularities of $f$ and the Kontsevich change of variables in motivic integration) that the real zeta function $Z_f^{?}$ is a rational function that can be expressed as

$$Z_f^?(T)= \sum_{I\cap \cal K \not= \emptyset}(\L-1)^{\vert I\vert -1 }[\widetilde E^{0,?}_I ]\prod_{i\in I}
\frac{\L^{-\nu_i}T^{N_i}}{1-\L^{-\nu_i}T^{N_i}} \eqno(12')$$
for $?$ being $-1,+1,>$ or $<$, where  $\widetilde E^{0,\epsilon}_I$ is defined as the gluing along $E^0_I\cap U$ of the sets  
$$ \{ (x,t)\in (E^0_I\cap U)\times \R; \ t^m \cdot u(x)\ !_? \ \},$$
where $!_?$ is $=-1$, $=1$, $>0$ or $<0$ in case $?$ is respectively $-1,+1, >$ or 
$< $.
 The real motivic Milnor $?$-fibre $S^?_f$ of $f$ may finally be defined as 
$$ S^?_f:=-\lim_{T\to \infty} Z^\epsilon_f(T):=\displaystyle
-\sum_{I\cap \cal K\not= \emptyset} (-1)^{\vert I \vert}
[{\widetilde E}_I^{0,?}](\L-1)^{\vert I \vert-1 }
\in K_0(\Var_\R)\otimes \Z[\frac{1}{2}].$$

\begin{remark}
The class $S^?_f$, although having an expression in terms of the data coming from a chosen resolution of $f$, does not depend of such a choice, 
since the definition of $Z^?_f$ as nothing to do with any choice of a resolution. 
\end{remark}

\begin{remark} There is  no {\sl a priori} obvious reason,  from the definition of $Z_f^?(T)$,
that  the constant term $S_f^?$ in the power series in $T^{-1}$ induced by the 
rational expression of $Z_f^?(T)$
could be accurately related to the topology of the corresponding set-theoretic Milnor fibre $X_0^?$, 
that is to say that $S_f^?$ could  be the motivic version 
of the signed Milnor fibre $X_0^?$ of $f$.  
In the complex case, it has just been observed that $\chi(S_f)$ 
is the expression of $\chi(X_0)$ provided by the A'Campo formula.  In the real case, taking into account 
that $\chi(\RR)=-1$, the 
expression of $\chi(S_f^?)$ is 
$$\chi(S_f^?)= \sum_{I\cap \cal K\not= \emptyset} (-2)^{\vert I \vert-1 }  \chi({\widetilde E}_I^{0,?}), $$
showing a greater complexity than in the complex case where only  strata $E_{\{i\}}$ of maximal dimension 
in the exceptional divisor $\sigma^{-1}(0)$ appear. 
Despite this increased complexity, in the real case the correspondence still holds, since it is 
proved in \cite{CoFi} that $\chi(S_f^?)$  is still $\chi(\bar X_0^?)$, $?\in \{-1,+1\}$. This justifies the terminology of  
motivic real semialgebraic Milnor fibre of $f$ at $0$ for $S_f^?$, at least at the first topological level 
represented by the morphism $\chi: K_0(\Var_\R)\otimes \Z[\frac{1}{2}] \to \ZZ$. 
\end{remark}

In order to accurately  state the correspondence between the motivic real semialgebraic Milnor fibre and the 
set-theoretic Milnor fibre we set now the following notation.  

\begin{notation}\label{Khim} Let us denote
$Lk(f)$ the link $f^{-1}(0)\cap S(0,\eta)$ of $f$ at the origin, $0<\eta\ll 1$. 
We recall that the topology of $Lk(f)$ is the same as the topology of the boundary 
$f^{-1}(\varepsilon)\cap S(0,\eta)$, $0<\varepsilon\ll \eta$,  of the Milnor fibre  $ f^{-1}(\varepsilon)\cap B_{(0,\eta)}$,
when $f$ has an isolated singularity at $0$.

- Let us denote, for $?\in \{<,>\}$, the topological type
of $f^{-1}(]0,c_?[)\cap B(0,\eta)$ by $X_0^?$, and 
the topological type of $f^{-1}(]0,c_?[)\cap \bar B(0,\eta)$
by $\bar X_0^?$,
where $c_<\in ]-\eta, 0[$ and $c_>\in ]0,\eta[$.

- Let us denote, for $?\in \{<,>\}$, the topological type
of $\{f \ \bar ? \ 0\}\cap S(0,\eta)$ by $G_0^?$, where 
$\bar ?$ is $\le$ when  $?$ is $<$ and $\bar ?$ is $\ge$
when  $?$ is $>$. 
\end{notation}

\begin{remark}\label{bord}
When $n$ is odd, $Lk(f)$ is a smooth odd-dimensional submanifold of $\RR^n$ and 
consequently $\chi(Lk(f))=0$. For 
$? \in \{-1, +1, <,>\}$, we thus have in this situation,  
$\chi(X_0^?)=\chi(\bar X_0^?)$. This is the situation in the complex 
setting. When $n$ is even, since
$\bar X_0^?$ is a compact manifold with boundary
$Lk(f)$, one knows from general algebraic topology that
$$\chi(\bar X_0^?)=-\chi (X_0^?)=\frac{1}{2}\chi(Lk(f)),$$
 for  $? \in \{-1, +1,<,>\}$.
For general $n\in \N$ and for $? \in \{-1, +1, <, >\}$, we thus have 
$$ \chi(\bar X_0^?)=(-1)^{n+1}\chi (X_0^?).$$
On the other hand we recall that for $? \in \{ <, >\}$
$$ \chi(G_0^?)=\chi(\bar X_0^{\delta_?}),$$
where $\delta_>$ is $+$ and $\delta_<$ is $-$ 
(see \cite{Ar}, \cite{Wa}).
\end{remark}

We can now state the real version of $(13)$.
We have, for $?\in \{-1, +1, <, >\}$ 
$$ \chi(S_f^?)=\sum_{I\cap \cal K\not= \emptyset} (-2)^{\vert I \vert-1 }  \chi({\widetilde E}_I^{0,?})=\chi(\bar X_0^?)=(-1)^{n+1}\chi (X_0^?),  \eqno(13')$$
and for $? \in \{<, >\}$
$$ \chi(S_f^?)=\sum_{I\cap \cal K\not= \emptyset} (-2)^{\vert I \vert-1 }  \chi({\widetilde E}_I^{0,?})=-\chi(G_0^?). \eqno(13'')$$

The formula 
$(13')$  below 
is the real analogue of the A'Campo-Denef-Loeser formula $(13)$ for complex hypersurface 
singularities and thus appears as the extension to the 
reals of this complex formula, or, in other words, the complex formula is the notably 
first level of complexity of the more general real formula $(13')$.

\begin{remark}
In \cite{Yin}, following the construction of Hrushovski and Kazhdan (see \cite{HruKaz1}, \cite{HruKaz2}), Yin develops a theory of motivic integration for polynomial bounded $T$-convex valued fields and studies, in this setting, topological 
zeta functions attached to a function germ, showing that they are rational. This a first step 
towards a real version of Hrushovski and Loeser work \cite{HruLoe}, where no resolution of singularities is used, in contrast with \cite{CoFi}. 
\end{remark}

\begin{questions}
The question of finding a real analogue of the complex monodromy with real analogues of the invariants
$\Lambda(M^m)$ is open. Similarly the question of defining a convenient zeta function with coefficients in an adapted
Grothendieck ring in order to let appear invariants of type $\Lambda_i^{\ell oc}$ or $\sigma_j$ ($e(\cal P^i)$ in the complex case) from a
 rational expression of this zeta function also naturally arises. 
\end{questions}

\bibliographystyle{amsplain}

\begin{thebibliography}{SGA}

 



\bibitem{ACA} N. A'Campo,
     \textit{La fonction z\^eta d'une monodromie}, Comment. Math. Helvetici 50 (1975), 223-248
     
     
\bibitem{Ale1} S. Alesker,  \textit{ On P. McMullen's conjecture on translation invariant valuations},
Adv. Math. 155 (2000), 239-263

\bibitem{Ale4} S. Alesker,  \textit{The multiplicative structure on continuous polynomial valuations},
GAFA, Geom. Funct. Anal. 14 (1) (2004), 1-26


\bibitem{Ale5} S. Alesker,  \textit{Theory of valuations on manifolds: a survey}, GAFA, Geom.
Funct. Anal. 17 (2007), 1321-1341

\bibitem{Ale-Fu} S. Alesker, J. H. G. Fu,  \textit{Theory of valuations on manifolds, III.
Multiplicative structure in the general case}, Trans. Amer. Mathematical Soc. 360 (4)
(2008), 1951-1981

\bibitem{Ar} V.I. Arnol'd,
   \textit{Index of a singular point of a vector fields, the Petrovski-Oleinik inequality, and
mixed Hodge structures}, Funct. Anal. and its Appl. 12 (1978), 1-14 
   
 \bibitem{Bat} V. Batyrev,  \textit{Birational Calabi-Yau n-folds have equal Betti numbers},
New trends in algebraic geometry (Warwick, 1996), London Math. Soc. Lecture Note Ser.,
264, Cambridge Univ.Press, Cambridge, (1999), 1-11
   
\bibitem{BeBr1} A. Bernig, L. Br\"ocker, \textit{Lipschitz-Killing invariants}, Math. Nachr. 245
(2002), 5-25
\bibitem{BeBr2} A. Bernig, L. Br\"ocker, \textit{Courbures intrins\`eques dans les cat\'egories
analytico-g\'eom\'etriques}, Ann. Inst. Fourier 53 (2003), no. 6, 1897-1924

\bibitem{Bli} M. Bickle, \textit{A short course on geometric motivic integration},
  Motivic integration and its interactions with model theory and non-Archimedean geometry, Volume I, 
London Math. Soc. Lecture Note Ser., 383, Cambridge Univ. Press, Cambridge, Edited by R. Cluckers, J. Nicaise, J. Sebag, 
(2011), 189-243

\bibitem{BrKu} L. Br\"ocker, M. Kuppe, \textit{Integral geometry of tame set}, Geom. Dedicata
82 (2000), no. 1-3, 285-323

\bibitem{BrDuKa} J. L. Brylinski, A. S. Dubson, M. Kashiwara, \textit{Formule de l'indice pour
modules holonomes et obstruction d'Euler locale}, C. R. Acad. Sci. Paris S\'er. I Math.
293 (1981), no. 12, 573-576

\bibitem{Car} P. Cartier, \textit{D\'ecomposition des poly\`edres : le point sur le troisi\`eme probl\`eme de Hilbert}, 
S\'eminaire Bourbaki, 1984-1985, exp. no 646, p. 261-288

\bibitem{Cha} Z. Chatzidakis, \textit{Introduction notes on the model theory of valued fields},
  Motivic integration and its interactions with model theory and non-Archimedean geometry, Volume I, 
London Math. Soc. Lecture Note Ser., 383, Cambridge Univ. Press, Cambridge, Edited by R. Cluckers, J. Nicaise, J. Sebag, 
(2011), 189-243

\bibitem{CheLas} S. S. Chern, R. Lashof, \textit{On the total curvature of immersed manifolds}, Amer. J.
Math., 79 (1957), 306-318


\bibitem{Clu-Loe2008} R. Cluckers, F. Loeser, \textit{Constructible motivic functions and motivic integration} Invent. Math. 173 (2008), no 1, 23-121

\bibitem{Clu-Loe2011} R. Cluckers, F. Loeser, \textit{Motivic integration in mixed characteristic with bounded ramification: a summary} Motivic integration and its interactions with model theory and non-Archimedean geometry, Volume I, 305-334, London Math. Soc. Lecture Note Ser., 383, Cambridge Univ. Press, Cambridge, (2011)



\bibitem{CoCRAS} G. Comte, \textit{Formule de Cauchy-Crofton pour la 
densit\'e des ensembles sous-analytiques}
C. R. Acad. Sci. Paris, t. 328 (1999), S\'erie I, 505-508

\bibitem{Co} G. Comte, \textit{\'Equisingularit\'e r\'eelle : nombres de Lelong et images
polaires},  Ann. Scient. \'Ec. Norm. Sup. 33(6) (2000), 757-788

\bibitem{CoFi} G. Comte, G. Fichou, 
\textit{Grothendieck ring of semialgebraic formulas and motivic real Milnor fibres}, (2013)
arXiv:1111.3181

\bibitem{CoMe} G. Comte, M. Merle, \textit{\'Equisingularit\'e r\'eelle II : invariants locaux et
conditions de r\'egularit\'e},   Ann. Scient. \'Ec. Norm. Sup. 41(2) (2008), 757-788
     
 \bibitem{CoLiRo}    G. Comte, J. -M. Lion, J.-Ph. Rollin, \textit{Nature Log-analytique du
volume des sous-analytiques}, Illinois J. Math 44, (4) (2000), 884-888

 \bibitem{De}   P. Deligne, \textit{Le formalisme des cycles \'evanescents}, S\'eminaire de G\'eom\'etrie Alg\'ebrique du Bois Marie, SGA7  XIII, 1967-69, Lecture Notes in Mathematics 340 (1973)
 
  \bibitem{Den84}   J. Denef, \textit{The rationality of the Poincar\'e series associated to the $p$-adic points on a variety}, Invent. Math. 77 (1984), no. 1, 1-23   
 
  \bibitem{DenBour}   J. Denef, \textit{Report on Igusa's local zeta function} 
S\'eminaire Bourbaki, Vol. 1990/1991, 
Ast\'erisque no. 201-203 (1991), Exp. no. 741, 359-386 (1992)


 \bibitem{DenLoe1992}   J. Denef, F. Loeser, \textit{Caract\'eristiques d'Euler-Poincar\'e, fonctions z\^eta locales et modifications analytiques}, J. Amer. Math. Soc. 5 (1992), no. 4, 705-720
 
  \bibitem{DenLoe1998}   J. Denef, F. Loeser, \textit{Motivic Igusa zeta functions} 
  J. Algebraic Geom. 7 (1998), no. 3, 505-537
 
 \bibitem{DenLoe1999}   J. Denef, F. Loeser, \textit{Germs of arcs on singular algebraic varieties and motivic integration} Invent. Math. 135 (1999), no. 1, 201-232

 \bibitem{DenLoe2001a}   J. Denef, F. Loeser, \textit{Definable sets, motives and p-adic integrals} J. Amer. Math. Soc. 14 (2001), no. 2, 429-469

  \bibitem{DenLoe2001b}   J. Denef, F. Loeser, \textit{Geometry on arc spaces of algebraic varieties}, European Congress of Mathematics, Vol. I (Barcelona, 2000), 327-348, Progr. Math., 201, Birkh\"auser, Basel, (2001)

  \bibitem{DenLoe2002}   J. Denef, F. Loeser, \textit{Lefschetz numbers of iterates of the monodromy and truncated arcs}, Topology 41 (2002), no. 5, 1031-1040

 \bibitem{Dra} R. N. Draper, \textit{Intersection theory in analytic geometry} Math. Ann. 180 (1969), 175-204
 
 
 \bibitem{Dub1}  A. S. Dubson, \textit{Classes caract\'eristiques des 
 vari\'et\'es singuli\`eres} C. R. Acad.
Sci. Paris S\'er. A-B 287 (1978), no. 4, 237-240

 \bibitem{Dub2} A. S. Dubson, \textit{Calcul des invariants num\'eriques des singularit\'es et des
applications} Th\`ese, Bonn University, (1981)

\bibitem{Du02} N. Dutertre, \textit{Courbures et singularit\'es r\'eelles}, Comment. Math. Helv. 77(4)
(2002), 846-863

\bibitem{Du08} N. Dutertre, \textit{A Gauss-Bonnet formula for closed semi-algebraic sets}, Advances in Geometry 8, no 1 (2008), 33-51

\bibitem{Du08Comm} N. Dutertre,  \textit{Curvature integrals on the real Milnor fiber}, Comment. Math.
Helvetici 83 (2008), 241-288

\bibitem{Du12} N. Dutertre, \textit{Euler characteristic and Lipschitz-Killing curvatures of closed semi-algebraic sets}, Geom. Dedicata 158 (2012), 167-189

\bibitem{Dri} L. van den Dries, \textit{Tame topology and o-minimal structures}, London 
Mathematical Society Lecture Note Series, 248, Cambridge University Press, Cambridge, (1998)
     
\bibitem{Fed1} H. Federer, \textit{ The $(\Phi, k)$ rectifiable subsets of n space}, Trans. Amer. Math.
Soc. 62 (1947), 114-192
     
     
\bibitem{Fed3} H. Federer,  \textit{Geometric measure theory}, Grundlehren Math. Wiss., Vol. 153
(1969) Springer-Verlag     

\bibitem{Fen} W. Fenchel, \textit{On total curvature of riemannian manifolds I},
Journal of London Math. Soc. 15 (1940), 15-22

\bibitem{Fi} G. Fichou, \textit{Motivic invariants of Arc-Symmetric sets and Blow-Nash Equivalence},
Compositio Math. 141 (2005), 655-688


\bibitem{FuSco} Joseph H.G. Fu, Ryan C. Scott, \textit{Piecewise linear approximation of smooth functions of two variables}, arXiv:1305.2220



\bibitem{Fu1} H. G. J. Fu, \textit{Tubular neighborhoods in Euclidean spaces}, Duke Math. J. 52
(1985), no. 4, 1025-1046

\bibitem{Fu2} H. G. J. Fu, \textit{Curvature measures and generalized Morse theory}, J. Differential
Geom. 30 (1989), no. 3, 619-642

\bibitem{Fu3}  H. G. J. Fu,  \textit{Monge-Amp`ere functions I}, Indiana Univ. Math. J. 38 (1989), 745-771

\bibitem{Fu4} H. G. J. Fu, \textit{ Monge-Amp`ere functions II}, Indiana Univ. Math. J. 38 (1989),
773-789

\bibitem{Fu5}  H. G. J. Fu, \textit{Kinematic formulas in integral geometry}, Indiana Univ. Math.
J. 39 (1990), no. 4, 1115-1154

\bibitem{Fu6}  H. G. J. Fu, \textit{Curvature of singular spaces via the normal cycle}, Differential
geometry: geometry in mathematical physics and related topics (Los Angeles, CA,
1990), 211-221, Proc. Sympos. Pure Math., 54 (1993), Part 2, Amer. Math. Soc.,
Providence, RI

\bibitem{Fu7} H. G. J. Fu, \textit{Curvature measures of subanalytic sets}, Amer. J. Math. 116
(1994), no. 4, 819-880

\bibitem{GorYaf} J. Gordon, Y. Yaffe, \textit{An overview of arithmetic motivic integration}
 Ottawa lectures on admissible representations of reductive p-adic groups, 113-149, 
Fields Inst. Monogr., 26, Amer. Math. Soc., Providence, RI, (2009)

\bibitem{Gr-Sc} P. M. Gruber, R. Schneider, \textit{Problems in geometric convexity. In:
Contributions to Geometry}, ed. par J. T\"olke et J. M. Wills, Birkh\"auser Verlag, Basel,
(1979), 225-278

\bibitem{GSLuMe} S. M. Gusein-Zade, I, Luengo,  A. Melle-Hern\'andez, 
\textit{Integration over a space of non-parametrized arcs, and motivic analogues of the monodromy zeta function}  Tr. Mat. Inst. Steklova 252 (2006), Geom. Topol., Diskret. Geom. i Teor. Mnozh., 71-82,
 translation in Proc. Steklov Inst. Math. 2006, no 1 (252), 63-73

\bibitem{Had} H. Hadwiger, \textit{Vorlesungen \"uber Inhalt, Oberfl\"ache und Isoperimetrie},
Springer-Verlag, Berlin-G\"ottingen-Heidelberg (1957)

\bibitem{Hal1} T. Hales, \textit{Can $p$-adic integrals be computed?} Contributions to automorphic forms, geometry, and number theory, 313-329, Johns Hopkins Univ. Press, Baltimore, MD, (2004)

\bibitem{Hal2} T. Hales, \textit{What is motivic measure?}  Bull. Amer. Math. Soc. (N.S.) 42 (2005), no. 2, 119-135

\bibitem{HeMe1} J. P. Henry, M. Merle, \textit{Limites de normales, conditions de Whitney et
\'eclatement d'Hironaka}, Proc. Symp. in Pure Math. 40 (1983) (vol. 1), Arcata 1981,
Amer. Math. Soc., 575-584

\bibitem{HruKaz1} E. Hrushovski, D. Kazhdan, \textit{Integration in valued fields, in Algebraic geometry
and number theory}, Progress in Mathematics 253, Birkh\"auser, (2006), 261-405

\bibitem{HruKaz2}  E. Hrushovski, D. Kazhdan, \textit{The value ring of geometric motivic integration, and
the Iwahori Hecke algebra of SL2. With an appendix by Nir Avni}, Geom. Funct.
Anal. 17 (2008), 1924-1967

\bibitem{HruLoe} E. Hrushovski, F. Loeser, \textit{Monodromy and the Lefschetz fixed point formula}
arXiv:1111.1954 


\bibitem{Igu1}  J. Igusa, \textit{Forms of higher degree} Tata Institute of Fundamental Research Lectures on Mathematics and Physics, 59, Tata Institute of Fundamental Research, Bombay, Narosa Publishing House, New Delhi, (1978)

\bibitem{Igu2}  J. Igusa, \textit{Complex powers and asymptotic expansions} II. Asymptotic expansions. J. Reine Angew. Math. 278/279 (1975), 307-321

\bibitem{Igu3}  J. Igusa, \textit{An introduction to the theory of local zeta functions} 
AMS/IP Studies in Advanced Mathematics, 14. American Mathematical Society, Providence, RI, International Press, Cambridge, MA, (2000)




\bibitem{Ka} M. Kashiwara, \textit{Index theorem for a maximally overdetermined system of
linear differential equations}, Proc. Japan Acad. 49 (1973), 803-804


\bibitem{Kla} D. A. Klain, \textit{A short proof of Hadwiger's characterization theorem},
Mathematika 42 (1995), 329-339

\bibitem{Kon} M. Kontsevich, \textit{Lecture at Orsay}, (December, 7 1995)

\bibitem{KonSoi} M. Kontsevich, Y. Soibelman \textit{Deformation Theory I}, preliminary draft 
http://www.math.ksu.edu/~soibel/Book-vol1.ps

\bibitem{KuRa} K. Kurdyka, G. Raby, \textit{ Densit\'e des ensembles sous-analytiques}, Ann. Inst.
Fourier 39 (1989), no. 3, 753-771
  
\bibitem{KuPoRa} K. Kurdyka, J. -P. Poly, G. Raby, \textit{Moyennes des fonctions sousanalytiques,
densit\'e, c\^one tangent et tranches}, (Trento, 1988), 170-177, Lecture Notes in Math., 1420 (1990), Springer, Berlin


\bibitem{Lan} R. Langevin, \textit{Courbure et singularit\'es complexes}, Comment.
Math. Helvetici 54 (1979), 6-16

\bibitem{Lan2} R. Langevin, \textit{Singularit\'es complexes, points critiques et int\'egrales de courbure}, S\'eminaire P. Lelong-H. Skoda, 18\`eme-19\`eme ann\'ee, (1978-1979), 129-143
 
\bibitem{LanLe} R. Langevin, L\^e Dung Tr\`ang,
\textit{ Courbure au voisinage d'une singularit\'e},  
C. R. Acad. Sci. Paris S\'er. A-B 290 (1980), no. 2
     
\bibitem{LanShi} R. Langevin, Th. Shifrin, \textit{Polar varietes and integral geometry}, Amer. J.  
Math. 104 (1982), no 3, 553-605 


\bibitem{LeTei1} L\^e D\~ung Tr\'ang, B. Teissier, \textit{Vari\'et\'es polaires locales et classes de Chern
des vari\'et\'es singuli\`eres}, Annals of Math. 114 (1981), 457-491


\bibitem{LeTei2} L\^e D\~ung Tr\'ang, B. Teissier, \textit{Errata \`a Vari\'et\'es polaires locales et classes
de Chern des vari\'et\'es singuli\`eres}, Annals of Math. 115 (1982), 668-668

\bibitem{LeTei3} L\^e D\~ung Tr\'ang, B. Teissier, \textit{Cycles \'evanescents et conditions de Whitney},
II. Proc. Symp. in Pure Math. 40 (1983) (vol. 2), Arcata 1981, Amer. Math. Soc.,
65-103


\bibitem{Li} J.-M. Lion, \textit{Densit\'e des ensembles semi-pfaffiens}, Ann. Fac. Sci. Toulouse Math. 6, 7 (1998), 
no. 1, 87-92    
     
\bibitem{Loe84}  F. Loeser, \textit{Formules int\'egrales pour certains invariants locaux des espaces analytiques complexes}, Comment. Math. Helv.
59 (1984), no. 2, 204-225  

\bibitem{Loe2009}  F. Loeser, \textit{Seattle lectures on motivic integration} Algebraic geometry-Seattle 2005, Part 2, 745-784, 
Proc. Sympos. Pure Math., 80, Part 2, Amer. Math. Soc., Providence, RI, (2009)

 

\bibitem{Loo} E. Looijenga, \textit{Motivic measures} 
S\'eminaire Bourbaki, Vol. 1999/2000, 
Ast\'erisque no. 276 (2002), 267-297


\bibitem{Mik} G. Mikhalkin,   \textit{Decomposition into pairs-of-pants for complex algebraic hypersurfaces}, 
Topology 43 (2004), no 5, 1035-1065

\bibitem{McMu-Sch} P. McMullen, R. Schneider, \textit{Valuations on convex bodies, Convexity
and its applications}, edited by Peter Gruber and J\"org M. Wills, Boston: Birkh\"auser
Verlag (1983)

\bibitem{MacPh}  R. MacPherson, \textit{Chern classes for singular algebraic varieties}, Ann. of Math. (2) 100
(1974), 423-432

 \bibitem{Me} M. Merle, \textit{Vari\'et\'es poalires, stratifications de Whitney et classes de Chern des espaces analytiques complexes} (d'apr\`es L\^e-Teissier), S\'eminaire Bourbaki, Vol. 1982/83, Exp. no 600, 
 Ast\'erisque 105 (1983), 65-78
 
     
\bibitem{Mil} J. Milnor, \textit{Singular points of complex hypersurfaces}, Ann. of Math. Studies 61 (1968)


\bibitem{Na2} V. Navarro Aznar, \textit{Stratifications r\'eguli\`eres et vari\'et\'es polaires locales}
Manuscrit, (1981)

\bibitem{NicSeb} J. Nicaise, J. Sebag, \textit{The Grothendieck rings of varieties},
  Motivic integration and its interactions with model theory and non-Archimedean geometry, Volume I, 
London Math. Soc. Lecture Note Ser., 383, Cambridge Univ. Press, Cambridge, Edited by R. Cluckers, J. Nicaise, J. Sebag, (2011), 145-188

\bibitem{RaiCRAS} M. Raibaut \textit{Fibre de Milnor motivique \`a l'infini}, C. R. Math. Acad. Sci. Paris, 348(7-8) (2010), 419-422

\bibitem{RaiSMF} M. Raibaut \textit{Singularit\'es \`a l'infini et int\'egration motivique}, 
Bull. SMF, 140(1) (2012), 51-100  


\bibitem{Rul} H. Rullg\aa rd, \textit{Polynomial amoebas and convexity}, preprint, Stockholm University, (2001)
Manuscrit, (1981)

\bibitem{Sab} C. Sabbah, \textit{Monodromy at infinity and Fourier transform}, Publ. Res. Inst. Math. Sci.,
33(4) (1997), 643-685

\bibitem{Sai} M. Saito, \textit{Mixed Hodge modules} Publ. Res. Inst. Math. Sci., 26(2) (1990), 221-333

 \bibitem{San} L. A. Santalo, \textit{Integral geometry and geometric probability}, Encyclopedia
of Mathematics and its Applications Vol. 1 (1976), Addison-Wesley Publishing Co.,London-Amsterdam


 \bibitem{Sch1} R. Schneider, \textit{Curvatures measures of convex bodies}, Ann. Mat. Pura appl.
116 (1978), 101-134
 
 \bibitem{Sch2} R. Schneider, \textit{A uniqueness theorem for finitely additive invariant measures
on a compact homogeneous space}, Rendiconti del Circolo Matematico di Palermo,
XXX (1981), 341-344


 \bibitem{Sch3} R. Schneider, \textit{Convex bodies: The Brunn-Minkowski Theory}, Encyclopedia
of Mathematics and its Applications 44 (1993), Cambridge University Press

 
 \bibitem{Sch4} R. Schneider, \textit{Integral geometry - Measure theoretic approach and stochastic
applications Advanced course on integral geometry}, CRM (1999)

 
 \bibitem{SGA7} SGA 7, 
 \textit{S\'eminaire de G\'eom\'etrie Alg\'ebrique du 
 Bois-Marie 1967-1969, Groupes de monodromie en g\'eom\'etrie alg\'ebrique (SGA 7), Vol. 1}, Springer Lecture Notes in Math. 288 (1972)
 
 
 \bibitem{SteZuc}  J. Steenbrink and S. Zucker, \textit{Variation of mixed Hodge structure. I}, Invent. Math.,
80(3) (1985), 489-542 

 \bibitem{St}  J. Steiner, \textit{\"Uber parallele Fl\"achen},
Monatsber. Preu$\beta$ Akad. Wissen., Berlin, (1840),
Ges. Werke, vol 2 (1882), Reimer, Berlin

\bibitem{Tei73} B. Teissier,  \textit{Cycles \'evanescents, sections panes et conditions de Whitney}, Ast\'erisque 7-8, Soc. Math. France (1973), 285-362


 \bibitem{Tei75}  B. Teissier, \textit{Introduction to equisingularity problems}, Proc. AMS  Symp. in Pure Math. 29, Arcata 1974, (1975) 
 
 \bibitem{Tei2}  B. Teissier, \textit{Vari\'et\'es polaires II : Multiplicit\'es polaires, sections planes et
conditions de Whitney} Actes de la conf\'erence de g\'eom\'etrie alg\'ebrique de la R\'abida
(1981), Springer Lecture Notes in Math. 961, Springer, Berlin, (1982), 314-491


 \bibitem{Tro1}  D. Trotman, \textit{Lectures on real stratification theory}, Singularity theory, World Sci. Publ., Hackensack, NJ, (2007), 139-155

 \bibitem{Tro2}  D. Trotman, \textit{Espaces stratifi\'es r\'eels}, Stratifications, singularities and differential equations Vol. II 
 (Marseille, 1990; Honolulu, HI, 1990),  Travaux en Cours 55, Hermann, Paris, (1997), 93-107


\bibitem{Veys}  W. Veys, \textit{Arc spaces, motivic integration and stringy invariants},   Singularity theory and its applications, 529-572, Adv. Stud. Pure Math., 43, Math. Soc. Japan, Tokyo, 2006

 \bibitem{Wa} C.T.C. Wall, \textit{Topological invariance of the Milnor number mod 2},
     Topology 22 (1983), 345-350 

  \bibitem{Wei}  A. Weil, \textit{Sur la formule de Siegel dans la th\'eorie des groupes classiques}, Acta Math. 113 (1965), 1-87


\bibitem{Wey} H. Weyl, \textit{On the Volume of Tubes}, Amer. J. Math. 61 (1939)

\bibitem{Yin} Y. Yin, 
\textit{Additive invariants in o-minimal valued fields}, (2013)
arXiv:1307.0224


\end{thebibliography}

\end{document}